\documentclass[12pt]{article}
\usepackage{amsmath, amssymb, xcolor, hyperref}

\newlength{\abstractwidth}
\renewcommand{\thefootnote}{\fnsymbol{footnote}}
\renewcommand{\thanks}[1]{\footnote{#1}} 
\newcommand{\starttext}{ \setcounter{footnote}{0}
\renewcommand{\thefootnote}{\arabic{footnote}}}

\newcommand{\be}{\begin{equation}}
\newcommand{\bea}{\begin{eqnarray}}
\newcommand{\eea}{\end{eqnarray}} \newcommand{\ee}{\end{equation}}

\def\ba{\begin{eqnarray}}
\def\ea{\end{eqnarray}}

\def\o{\omega}

\def\det{{\rm det}}

\def\log{\,{\rm log}\,}

\def\o{\omega}

\def\o{\omega}

\def\na{\nabla}

\def\p{\prod}

\def\na{{\nabla}}

\def\[{{\bf [}}
\def\]{{\bf ]}}

\def\p{\partial}

\begin{document}
\starttext \baselineskip=15pt \setcounter{footnote}{0}
\newtheorem{theorem}{Theorem}
\newtheorem{lemma}{Lemma}
\newtheorem{definition}{Definition}
\newtheorem{proposition}{Proposition}

\begin{center}
{\Large \bf Bochner-Kodaira Formulas and the Type IIA Flow
\footnote{Work supported in part by the National Science Foundation Grants DMS-1855947 and DMS-1809582.}}
\bigskip\bigskip

\centerline{Teng Fei, Duong H. Phong, Sebastien Picard,
and Xiangwen Zhang}

\begin{abstract}
{\footnotesize A new derivation of the flow of metrics in the Type IIA flow is given. It is adapted to the formulation of the flow as a variant of a Laplacian flow, and it uses the projected Levi-Civita connection of the metrics themselves instead of their conformal rescalings.}

\end{abstract}

\end{center}

\section{Introduction}
\setcounter{equation}{0}

The search for supersymmetric compactifications of string theories has revealed itself to have deep connections with special geometry. The resulting non-linear partial differential equations also turned out to be quite rich and interesting in their own right
(see e.g. \cite{FHP, FGV, FY1,FY2,PPZ0, PPZ3}). One invariable feature of particular interest in these equations is the presence of a cohomological constraint. When no $\p\bar\p$-lemma is available, the most natural implementation of these cohomological constraints is by a geometric flow, and this has resulted in considerable interest in the investigation of such flows in recent years \cite{BV, B,BX, FPPZIIB, FFM, Lo, LoW,Ph, PPZ1,PPZ2}.

\medskip
The present paper is mainly concerned with the Type IIA flow, which is a flow in symplectic geometry  introduced in \cite{FPPZIIA} and motivated by the Type IIA string. More specifically, let $(M,\o)$ be a compact $6$-dimensional symplectic manifold and $\rho_A$ be the Poincar\'e  dual to a finite combination of Lagrangians. Then the Type IIA flow is the flow of $3$-forms $\varphi$ given by
\bea
\label{eq:IIA}
\p_t\varphi=d\Lambda d(|\varphi|^2\star\varphi)-\rho_A
\eea
with an initial data $\varphi_0$ which is a closed, primitive, and positive $3$-form on $M$. Here $\Lambda$ is the Hodge contraction operator defined by $\o$, and $\star$ and $|\varphi|$ are the Hodge star operator and the norm of $\varphi$ with respect to the metric $g_\varphi$ which is compatible with $\o$ and the almost-complex structure $J_\varphi$ constructed by Hitchin \cite{Hitchin} (see \S \ref{background} for the precise definitions). The Type IIA flow preserves the primitiveness and closedness of $\varphi$, so that its stationary points are automatically solutions of the system first investigated by Tseng and Yau \cite{TY3}. This system is itself a basic case of the more general equations for supersymmetric compactifications of the Type IIA string proposed in \cite{Getal,T}.

\medskip
In \cite{FPPZIIA} it was shown that the Type IIA flow admits at least short-time existence, and can be continued as long as $|\varphi|$ and the Riemannian curvature of $g_\varphi$ remain bounded. The proof of this last assertion relied heavily on determining the flow of $g_\varphi$. This was one of the main results of \cite{FPPZIIA}, and it was established using the original formulation (\ref{eq:IIA}) of the Type IIA flow, and the projected Levi-Civita connection $\tilde{\frak D}$ of a metric $\tilde g_\varphi$ conformal to $g_\varphi$ (see (\ref{tildeg} below). A key point was that, with respect to $\tilde{\frak D}$, the manifold $M$ has SU(3) holonomy, and
the form $|\Omega|_{\tilde g_\varphi}^{-1}\Omega$, with $\Omega=\varphi+i\star\varphi$, is covariant constant.

\medskip
The main goal of the present paper is to provide a different derivation of the flow of the metrics $g_\varphi$ in the Type IIA flow. The new derivation differs from the one in \cite{FPPZIIA} in two important aspects. The first aspect is that it relies on Bochner-Kodaira formulas and a different formulation of the Type IIA flow, which is closer in spirit to Bryant's $G_2$ flow. From this point of view, it is more easily adaptable to other Laplacian flows.
The second aspect is that it relies instead on the projected Levi-Civita connection ${\frak D}$ of $g_\varphi$, which is a very natural connection since it coincides with all the unitary Hermitian connections with respect to $g_\varphi$ on the Gauduchon line. An important additional benefit of this second derivation is that it provides a check on the formulas obtained in \cite{FPPZIIA}, which is non-trivial because the calculations in both approaches are particularly long and involved.

\medskip
For simplicity, we focus on the source-free case $\rho_A=0$. Then we have

\begin{theorem}
\label{main}
Let $(M,\o)$ be a $6$-dimensional symplectic manifold, and let $t\to\varphi(t)$ by the Type IIA flow of $3$-forms defined in (\ref{eq:IIA}) with $\rho_A=0$. If $g_{ij}=(g_\varphi)_{ij}$ is the corresponding flow of metrics, then we have
\bea \label{evol-g}
\partial_t g_{ij} = - | \varphi|^2 \left\{ 2 R_{ij}   - 2  \nabla_i \nabla_j \log | \varphi |^2 + 4 (N^2_-)_{ij}  - \alpha_i \alpha_j  + \alpha_{Ji} \alpha_{Jj}   +4 \alpha_p (N_j{}^p{}_i + N_i{}^p{}_j) \right\}
\eea
where $\nabla$ is the Levi-Civita connection of $g$, $R_{ij}$ is the Ricci curvature, $N$ is the Nijenhuis tensor with respect to the almost-complex structure $J_{\varphi}$,
$(N_-^2)_{ij}=N^{\lambda p}{}_iN_{p\lambda j}$, and $\alpha$ is the 1-form defined by $\alpha = -d\log |\varphi|^2$.
\end{theorem}

\section{Background material}
\setcounter{equation}{0}
\label{background}

We begin by providing a brief summary of the setting for the Type IIA flow,
which is Type IIA geometry as introduced in \cite{FPPZIIA}.

\subsection{Type IIA geometry}

Let $M$ be an oriented $6$-manifold. In \cite{Hitchin}, Hitchin has shown how to associate to any non-degenerate $3$-form $\varphi$ an almost-complex structure $J_\varphi$. Type IIA geometry arises if, in addition, $M$ is equipped with a fixed symplectic form $\o$ and $\varphi$ is a closed form which is primitive and positive with respect to $\o$. The primitive condition means that $\Lambda\varphi=0$, where $\Lambda:A^k(M)\to A^{k-2}(M)$ is the standard Hodge contraction operator with respect to $\o$. It is shown in \cite{FPPZIIA} that $\o$ is then preserved by $J_\varphi$, and the positivity condition means that the resulting Hermitian form $g_\varphi(X,Y)=\o(X,J_\varphi Y)$ is positive definite and defines a metric. Thus $(J_\varphi,g_\varphi,\o)$ is an almost-K\"ahler manifold. However, the condition in Type IIA geometry that this almost-K\"ahler structure arise from a closed $3$-form results in many subtle properties which are essential for the Type IIA flow.

\medskip
Explicitly, the metric $g_\varphi$ is given by
\bea
(g_\varphi)_{ij}=-|\varphi|^{-2}\varphi_{iab}\varphi_{jkp}\o^{ak}\o^{bp}
\eea
where $|\varphi|$ is the norm of the $3$-form $\varphi$ with respect to $J_\varphi$, and $\o^{ak}$ is the inverse of the symplectic form $\o$, $\o^{ak}\o_{kp}=\delta^a{}_p$. The volume form of $g_\varphi$ is the same as $\o^3/3!$.
The following metric $\tilde g_\varphi$ conformally equivalent to $g_\varphi$ also plays an important role in Type IIA geometry,
\bea
\label{tildeg}
(\tilde g_\varphi)_{ij}=|\varphi|^2(g_\varphi)_{ij}=-\varphi_{iab}\varphi_{jkp}\o^{ak}\o^{bp}.
\eea
In fact, one of the defining features of Type IIA geometry is that the manifold $(M,J_\varphi)$ have SU(3) holonomy with respect to the projected Levi-Civita connection $\tilde{\frak D}$ of $\tilde g_\varphi$. More precisely, 
set
\bea
\hat\varphi=\star\varphi=J\varphi
\eea
and let $\Omega$ be the $(3,0)$-form defined by
\bea
\Omega=\varphi+i\hat\varphi.
\eea
Then $|\Omega|_{\tilde g_\varphi}^{-1}\,\Omega$ is covariantly constant with respect to $\tilde{\frak D}$. This was a major reason why the calculations in \cite{FPPZIIA} were mostly carried out with the connection $\tilde{\frak D}$.

\medskip
In the present paper, we shall use instead the unitary connections with respect to $g_\varphi$. Since $\o$ is closed, the Gauduchon line of Hermitian unitary connections with respect to $J_\varphi$ collapses to a single connection, which can be viewed as either the Chern connection or the projected Levi-Civita connection ${\frak D}$ of $g_\varphi$. Henceforth we drop the subindex $\varphi$ when there is no possibility of confusion, and denote $g_\varphi$, $\tilde g_\varphi$, $J_\varphi$ simply by $g$, $\tilde g$ and $J$. Then the Levi-Civita connection $\na$ and the projected Levi-Civita connection ${\frak D}$ of $g$ are related by
\be
\label{frakD-nabla}
{\frak D}_i X^m = \nabla_i X^m - N_{ip}{}^m X^p
\ee
where $N_{ip}{}^m$ is the Nijenhuis tensor of $J$,
\be
N^k{}_{ij} = {1 \over 4} (J^r{}_i \nabla_r J^k{}_j + J^k{}_r \nabla_j
J^r{}_i - (i \leftrightarrow j) ).
\ee
In \cite{FPPZIIA}, we showed ${\frak D}^{0,1}\Omega=0$ and ${\frak D}^{1, 0} \Omega = - \alpha\otimes \Omega$ (Equation (6.50) in \cite{FPPZIIA}), or equivalently,
\bea
\label{frakD-varphi}
{\frak D}_m\varphi={1\over 2}(-\alpha_m\varphi-\alpha_{Jm}\hat\varphi),
\quad
{\frak D}_m\hat\varphi=
{1\over 2}(-\alpha_m\hat\varphi+\alpha_{J m}\varphi).
\eea
Here the $1$-form $\alpha$ is defined by
\bea
\alpha=-d\log |\varphi|^2
\eea
and we used the same notation introduced in \cite{FPPZIIA} for any vector field $V$ and any 1-form $W$,
\be
V^{J k}=(JV)^k = J^k{}_p V^p, \quad W_{J k}=(JW)_k = J^p{}_k W_p.
\ee
In particular, $\omega_{ij} = g_{Ji,j}$,  $g_{ij} = \omega_{i,Jj}$, and $\omega^{ij} = g^{Ji,j}$,  $g^{ij} = \omega^{i,Jj}$.

\subsection{Identities from Type IIA geometry}

We list here some identities required later. Except for (\ref{riemann-J-identity}), they were proved in \cite{FPPZIIA}. 

\subsubsection{Identities for $\varphi$}

First, the action of $J$ on $\varphi$ is given by
\bea
&&
\varphi_{ijk}=-\varphi_{Ji,Jj,k}=-\varphi_{Ji,j,Jk}=-\varphi_{i,Jj,Jk}
\nonumber\\
&&
\varphi_{Ji,j,k}
=\varphi_{j,Jj,k}=\varphi_{i,j,Jk}.
\eea
Next, bilinears in $\varphi$ with two contractions with $\o^{ij}$ give the metric $g_{ij}$. But bilinears with a single contraction with either $\o^{ij}$ or $g^{ij}$ simplify as well,
\bea\label{contract}
&&
\o^{ij}\varphi_{iab}\varphi_{jcd}
=
{| \varphi |^2 \over 4}(\o_{ac}g_{bd}+\o_{bd}g_{ac}
-\o_{bc}g_{ad}-\o_{ad}g_{bc})
\nonumber\\
&&
g^{ij}\varphi_{iab}\varphi_{jcd}
=
{| \varphi |^2  \over 4}(g_{ac}g_{bd}+\o_{ca}\o_{bd}-\o_{ad}\o_{cb}-g_{bc}g_{ad}).
\eea
As a consequence, we also have bilinear identities involving $\varphi$ and $\hat\varphi$, for example
\bea
\hat\varphi_{\lambda kp}\varphi_{iab}\o^{ka}\o^{pb}
=
| \varphi |^2 \o_{\lambda i}.
\eea
This reduces to the previous identity by noting that $\hat\varphi_{\lambda kp}=-\varphi_{J\lambda,kp}$, so that
\bea
\hat\varphi_{\lambda kp}\varphi_{iab}\o^{ka}\o^{pb}=-\varphi_{J\lambda,kp}\varphi_{iab}\o^{ka}\o^{pb}
=| \varphi |^2 g_{J\lambda, i}=| \varphi |^2 \o_{\lambda i}.
\eea

\subsubsection{Identities for the Nijenhuis tensor}

In general, the Nijenhuis tensor satisfies the following identities of a type $(0,2)$-tensor in the sense of Gauduchon \cite{G}
\be
N^k{}_{Ji,j}= - N^{Jk}{}_{ij} = N^k{}_{i,Jj}, \quad N_{Ji,j,k} = N_{i,Jj,k} = N_{i,j,Jk}.
\ee
Since $d \omega = 0$, we also have the Bianchi identity
\be
N_{ijk} + N_{jki} + N_{kij} = 0.
\ee
From this it follows that there are two symmetric tensors quadratic in
$N$, denoted by
\be
(N^2_+)_{ij} = N^{pq}{}_i N_{pq j}, \quad (N^2_-)_{ij} = N^{pq}{}_i N_{qpj}.
\ee
The relation between the Levi-Civita connection $\na$ and the projected Levi-Civita connection ${\frak D}$ also implies, since ${\frak D}J=0$,
\be
\nabla_i J^k{}_j = - 2 N_{ij}{}^{Jk}.
\ee
In Type IIA geometry, we have
\bea\label{qua-N}
N_-^2=2\,N_+^2-{1\over 4}|N|^2 g,
\quad
|N|^2=(N_+^2)^\lambda{}_\lambda=2(N_-^2)^\lambda{}_\lambda
\eea
where $|N|^2=N^{mkp}N_{mkp}$. We also have the following crucial identity between the Nijenhuis tensor and $\varphi$,
\be \label{Nphi-switch}
N^p{}_{ij} \varphi_{pkl} = - N^p{}_{kl} \varphi_{p ij},
\ee
which was proved in \cite{FPPZIIA}, Corollary 1.

\subsubsection{Identities for the curvature tensor}

We shall express the desired identities for the curvature tensor of the Levi-Civita connection in the following convention. The connection $\na$ is written as
$\nabla_m V^k=\p_mV^k+\Gamma^k{}_{m\ell}V^\ell$, and the curvature tensor $R_{ij}{}^k{}_\ell$ is defined by
\bea
[\na_i,\na_j]V^k=R_{ij}{}^k{}_\ell V^\ell.
\eea
The Ricci curvature is then given by $R_{ij} = R_{ipj}{}^p$.

\medskip
The first curvature identity that we require gives the action of $J$ on $Rm$,
\bea \label{riemann-J-identity}
R_{j, i, Jk, J \ell} &=& R_{jik \ell}+B_{ijk\ell}
\nonumber\\
B_{ijk\ell}&=&- 2 {\frak D}_i N_{j k \ell } + 2 {\frak D}_j N_{i k \ell} - 2 N^\alpha{}_{ij} N_{\alpha k \ell} ,
\eea
This identity can also be expressed as 
\be
 R_{ji}{}^p{}_{J\ell} = R_{ji}{}^{Jp}{}_\ell + 2 {\frak D}_j  N_{i}{}^{Jp}{}_\ell - 2 {\frak D}_i  N_{j}{}^{Jp}{}_\ell -2  N^\mu{}_{ji} N_{\mu \ell}{}^{Jp}.
\ee
To see this, we consider the action of $J$ on a vector field $V$,
\bea
R_{jk}{}^p{}_q(JV)^q&=&\na_j\na_k(JV)^p-\na_k\na_j(JV)^p
\nonumber\\
&=&
J[\na_j,\na_k]V^p+(\na_j\na_kJ-\na_k\na_jJ)^p{}_\lambda V^\lambda.
\eea
It follows that
\bea
R_{jk}{}^p{}_qJ^q{}_\lambda
&=&
J^p{}_qR_{kj}{}^q{}_\lambda+\na_j\na_kJ^p{}_\lambda-\na_k\na_jJ^p{}_\lambda
\nonumber\\
&=&
J^p{}_qR_{kj}{}^q{}_\lambda-2\na_j (J^p{}_\mu  N_{k\lambda}{}^\mu ) +2\na_k ( J^p{}_\mu N_{j\lambda}{}^\mu)
\eea
or, in more succinct notation,
\be \label{riemman-J-identity0}
R_{jk}{}^p{}_{J\lambda}=R_{jk}{}^{Jp}{}_\lambda -2\na_j (N_{k\lambda}{}^{Jp}) +2\na_k (N_{j\lambda}{}^{Jp}).
\ee
We now convert $\nabla$ derivatives into ${\frak D}$
derivatives. First lowering indices gives
\be
R_{ji k,J \ell}= - R_{j,i,Jk,\ell} +2 \na_j (N_{i, \ell, Jk}) -2 \na_i (N_{j, \ell,Jk}).
\ee
Therefore
\be
R_{j,i, Jk,J \ell}=  R_{ji k \ell} +2 J^p{}_k \na_j (N_{i, \ell, Jp}) -2 J^p{}_k \na_i (N_{j, \ell,Jp}).
\ee
We write
\bea
2 J^p{}_k \na_j (N_{i, \ell, Jp}) &=& 2 J^p{}_k {\frak D}_j (N_{i, \ell, Jp}) -  2 J^p{}_k N_{ji}{}^\mu (N_{\mu, \ell, Jp}) \nonumber\\
&&-  2 J^p{}_k N_{j \ell}{}^\mu (N_{i, \mu, Jp}) - 2 J^p{}_k N_{jp}{}^\mu (J^n{}_\mu N_{i \ell n})
\eea
Since ${\frak D} J =0$,
\bea
2 J^p{}_k \na_j (N_{i, \ell, Jp}) &=& -2 {\frak D}_j  N_{i \ell k} +2  N_{ji}{}^\mu N_{\mu \ell k} + 2 N_{j \ell}{}^\mu N_{i \mu k} - 2 N_{j,Jk}{}^{Jn}  N_{i \ell n} \nonumber\\
 &=& 2 {\frak D}_j  N_{i k \ell} +2  N_{ji}{}^\mu N_{\mu \ell k} - 2( N_{j \ell}{}^\mu N_{i k \mu} + N_{jk}{}^{\mu}  N_{i \ell \mu})
\eea
This last term is symmetric in $(i,j)$. Therefore
\bea
2 J^p{}_k \na_j (N_{i, \ell, Jp}) - (i \leftrightarrow j)  
 &=& 2 {\frak D}_j  N_{i k \ell} - 2 {\frak D}_i  N_{j k \ell} +2  N_{ji}{}^\mu N_{\mu \ell k} - 2 N_{ij}{}^\mu N_{\mu \ell k}
\eea
By the Bianchi identity
\bea
2 J^p{}_k \na_j (N_{i, \ell, Jp}) - (i \leftrightarrow j)  
&=& 2 {\frak D}_j  N_{i k \ell} - 2 {\frak D}_i  N_{j k \ell} +2  (- N^\mu{}_{ji} - N_i{}^\mu{}_j) N_{\mu \ell k} \nonumber\\
&&- 2 N_{ij}{}^\mu N_{\mu \ell k}
\eea
from which the desired identity (\ref{riemann-J-identity}) follows.

\medskip
Finally, we shall need the following curvature identity specific to Type IIA geometry (see (6.53) in \cite{FPPZIIA}),
\be \label{ricci-typeiia}
R_{ij} = - {\frak D}_s (N_{i}{}^s{}_j + N_j{}^s{}_i)  - 2 (N^2_-)_{ij} + {1 \over 2} \nabla_i \nabla_j \log | \varphi |^2 + {1 \over 2} J^p{}_iJ^q{}_j \nabla_{p} \nabla_{q} \log | \varphi |^2.
\ee

\

\section{Proof of Theorem \ref{main}}
\setcounter{equation}{0}

We shall establish Theorem \ref{main} using the formulation of the Type IIA flow as a Laplacian type flow \cite{FPPZIIA}
\bea
\p_t\varphi= - dd^\dagger(|\varphi|^2\varphi) + 2 d ( | \varphi |^2 N^\dagger \cdot \varphi)
\eea
where $N^\dagger: \Lambda^3(M)\to\Lambda^2(M)$ is the operator defined by
\bea
\label{Ndagger}
(N^\dagger\cdot\varphi)_{kj}
=N^{\mu}{}_j{}^\lambda \varphi_{\mu k\lambda}-N^{\mu}{}_k{}^\lambda \varphi_{\mu j\lambda}.
\eea
For our present purposes, it is convenient to rewrite the above expression as
\bea \label{phi-flow-rhs}
\p_t\varphi= - |\varphi|^2 dd^\dagger\varphi -d|\varphi|^2\wedge d^\dagger\varphi + d(\iota_{\na|\varphi|^2}\varphi)
+2d(|\varphi|^2N^\dagger\cdot\varphi).
\eea
We would like to determine $\p_tg_{ij}$ explicitly. For this, it is convenient to determine first $\p_t\tilde g_{ij}$, since $\tilde g_{ij}$ is a quadratic expression in $\varphi$, and
we have
\bea\label{flow-tildeg}
\p_t\tilde g_{ij}=-\left\{(\p_t \varphi_{iab})\varphi_{jkp}\o^{ka}\o^{pb}+(i\leftrightarrow j)\right\}.
\eea
We shall determine in turn the contribution of each expression in (\ref{phi-flow-rhs}) to $\p_t\tilde g_{ij}$.

\subsection{The Bochner-Kodaira formula for the Levi-Civita connection}

We begin with the contribution of $|\varphi|^2dd^\dagger\varphi$ using a Bochner-Kodaira formula. In general, if $M$ is any compact Riemannian manifold and we express any $p$-form in components as
\bea
\varphi={1\over p!}\sum_{i_1,\cdots,i_p}\varphi_{i_1\cdots i_p}dx^{i_1}\wedge\cdots\wedge dx^{i_p}={1\over p!}\sum_{I}\varphi_{I}dx^I
\eea
with antisymmetric coefficients $\varphi_{i_1\cdots i_p}$, then the adjoint $d^\dagger$ of the de Rham exterior differential with respect to a given metric $g_{ij}$ is given by
\bea
(d^\dagger\varphi)_{I'}=-g^{\ell m}\nabla_m\varphi_{\ell I'},
\eea
where $\nabla$ denotes the covariant derivative with respect to the Levi-Civita connection and we have split the index $I$ into $I=(\ell,I')$,
$I'=(i_2,\cdots,i_p)$. It follows that
\bea
(dd^\dagger\varphi)_{I}
=
-(\na_{i_1}(g^{\ell m}\na_m\varphi_{\ell I'})
-
\sum_{q=2}^p(i_1\leftrightarrow i_q)).
\eea
Next, we have
\bea
(d\varphi)_{\ell I}
=
\na_\ell\varphi_{I}-\sum_{q=1}^p(\ell\leftrightarrow i_q)
\eea
and hence
\bea
(d^\dagger d\varphi)_I
=
-g^{\ell m}\na_m(\na_\ell\varphi_I-\sum_{q=1}^p(\ell\leftrightarrow i_q)).
\eea
Altogether, we obtain the version of the Bochner-Kodaira formula that we need,
\bea
((dd^\dagger+d^\dagger d)\varphi)_I
=
-g^{\ell m}\na_m\na_\ell\varphi_I
+
g^{\ell m}\sum_{q=1}^p[\na_m,\na_{i_q}]\varphi_{\cdots i_{q-1}\ell i_{q+1}\cdots}
\eea
In the case of interest, namely $3$-forms $\varphi$ with $d\varphi=0$, we obtain
\bea \label{bk-identity1}
dd^\dagger\varphi_{jkp}
=
-g^{\ell m}\na_m\na_\ell\varphi_{jkp}
+
g^{\ell m}
\big\{[\na_m,\na_j]\varphi_{kp\ell}
+
[\na_m,\na_k]\varphi_{pj \ell}
+
[\na_m,\na_p]\varphi_{jk\ell}\big\}.
\eea

\subsection{The Laplacian term $g^{\ell m}\na_m\na_\ell\varphi_{jkp}$}

Recall that the covariant derivatives of $\varphi$ with respect to the projected Levi-Civita connection ${\frak D}$ are given by (\ref{frakD-varphi}). It follows that
\bea
g^{\ell m}{\frak D}_\ell{\frak D}_m\varphi=- {1 \over 2} (\nabla_\mu \alpha^\mu)\varphi
\eea
and
\bea
[{\frak D}_m,{\frak D}_\ell]\varphi
&=&{1 \over 2} (-{\frak D}_m\alpha_\ell +{\frak D}_\ell \alpha_m)\varphi+ {1 \over 2} (-{\frak D}_m\alpha_{J \ell}+{\frak D}_\ell \alpha_{Jm})\hat\varphi
\nonumber\\
&=&
- {1 \over 2} \,N_{m\ell}{}^j\alpha_j\,\varphi + {1 \over 2} \,N_{\ell m}{}^j\alpha_j\,\varphi + {1 \over 2} (-{\frak D}_m\alpha_{J \ell}+{\frak D}_\ell \alpha_{Jm})\hat\varphi.
\eea
Now the difference between $\nabla$ and ${\frak D}$ on vectors is given by (\ref{frakD-nabla}). On $3$-forms, it is given by
\bea
\na_\ell\varphi_{jkp}&=&{\frak D}_\ell\varphi_{jkp}
-
\varphi_{\lambda kp}N_{\ell j}{}^\lambda
-
\varphi_{j\lambda p}N_{\ell k}{}^\lambda
-
\varphi_{jk\lambda}N_{\ell p}{}^\lambda
\nonumber\\
&=&
{\frak D}_\ell\varphi_{jkp}
-E_{\ell;jkp},
\eea
where
\be
E_{\ell;jkp} = \varphi_{\lambda kp}N_{\ell j}{}^\lambda + \varphi_{j\lambda p}N_{\ell k}{}^\lambda + \varphi_{jk\lambda}N_{\ell p}{}^\lambda.
\ee
Similarly, we write
\bea
\na_m{\frak D}_\ell\varphi_{jkp}&=&{\frak D}_m{\frak D}_\ell\varphi_{jkp}
-
{\frak D}_\mu\varphi_{jkp}N_{m\ell}{}^\mu
-
{\frak D}_\ell\varphi_{\mu kp}N_{mj}{}^\mu
-
{\frak D}_\ell\varphi_{j\mu p}N_{mk}{}^\mu
-
{\frak D}_\ell\varphi_{jk\mu}N_{mp}{}^\mu
\nonumber\\
&:=&
{\frak D}_m{\frak D}_\ell\varphi_{jkp}
-E_{m;\ell jkp},
\eea
and hence
\bea
g^{m\ell}\na_m\na_\ell\varphi_{jkp}
=
g^{m\ell}{\frak D}_m{\frak D}_\ell\varphi_{jkp}-g^{m\ell}E_{m;\ell jkp}-g^{m\ell}\na_m E_{\ell;jkp}.
\eea

\smallskip

We begin by computing the contributions of $g^{m\ell}\na_m E_{\ell;jkp}$,
\bea
\label{NaE}
g^{m\ell}\na_m E_{\ell;jkp}
&=&
(g^{\ell m}\na_m\varphi_{\lambda kp})N_{\ell j}{}^\lambda
+
(g^{\ell m}\na_m\varphi_{j\lambda p})N_{\ell k}{}^\lambda
+
(g^{\ell m}\na_m\varphi_{jk\lambda})N_{\ell p}{}^\lambda
\nonumber\\
&&
+
\varphi_{\lambda kp}g^{\ell m}\na_mN_{\ell j}{}^\lambda
+
\varphi_{j\lambda p}g^{\ell m}\na_mN_{\ell k}{}^\lambda
+
\varphi_{jk\lambda}g^{\ell m}\na_mN_{\ell p}{}^\lambda
\nonumber\\
&=&
(g^{\ell m}{\frak D}_m\varphi_{\lambda kp})N_{\ell j}{}^\lambda
+
(g^{\ell m}{\frak D}_m\varphi_{j\lambda p})N_{\ell k}{}^\lambda
+
(g^{\ell m}{\frak D}_m\varphi_{jk\lambda})N_{\ell p}{}^\lambda
\nonumber\\
&&
-
g^{\ell m}
(E_{m;\lambda kp}N_{\ell j}{}^\lambda
+
E_{m;j\lambda p} N_{\ell k}{}^\lambda
+
E_{m;jk\lambda}N_{\ell p}{}^\lambda)
\nonumber\\
&&
+
\varphi_{\lambda kp}g^{\ell m}\na_mN_{\ell j}{}^\lambda
+
\varphi_{j\lambda p}g^{\ell m}\na_mN_{\ell k}{}^\lambda
+
\varphi_{jk\lambda}g^{\ell m}\na_mN_{\ell p}{}^\lambda.
\eea

\subsubsection{Contributions of the terms $E_{\ell;jkp}$}

Consider the contributions of the second row on the right hand side of the last equation. Paired with $\varphi_{iab}\o^{ka}\o^{pb}$, it gives
\bea
g^{\ell m}E_{m;\lambda kp}N_{\ell j}{}^\lambda
\,\varphi_{iab}\o^{ka}\o^{pb}
&=&
g^{\ell m}(\varphi_{\mu kp}N_{m\lambda}{}^\mu
+
\varphi_{\lambda\mu p}N_{mk}{}^\mu
+
\varphi_{\lambda k\mu}N_{mp}{}^\mu)N_{\ell j}{}^\lambda
\,\varphi_{iab}\o^{ka}\o^{pb}
\nonumber\\
&=&
({\rm I}+{\rm II}+{\rm III})\cdot
\varphi_{iab}\o^{ka}\o^{pb}
\eea
with
\bea
({\rm I})\cdot\varphi_{iab}\o^{ka}\o^{pb}
&=&
g^{\ell m}\varphi_{\mu kp}N_{m\lambda}{}^\mu N_{\ell j}{}^\lambda\,\varphi_{iab}\o^{ka}\o^{pb}
\nonumber\\
({\rm II})\cdot\varphi_{iab}\o^{ka}\o^{pb}
&=&
g^{\ell m}\varphi_{\lambda\mu p}N_{mk}{}^\mu N_{\ell j}{}^\lambda\,\varphi_{iab}\o^{ka}\o^{pb}
\nonumber\\
({\rm III})\cdot\varphi_{iab}\o^{ka}\o^{pb}
&=&
g^{\ell m}\varphi_{\lambda k\mu}N_{mp}{}^\mu N_{\ell j}{}^\lambda\,\varphi_{iab}\o^{ka}\o^{pb}.
\eea
Next, we have
\bea
({\rm I})\cdot\varphi_{iab}\o^{ka}\o^{pb}
=
- | \varphi |^2 g^{\ell m}g_{\mu i}N_{m\lambda}{}^\mu N_{\ell j}{}^\lambda
=
- | \varphi |^2 N^\ell{}_{\lambda i}N_{\ell j}{}^\lambda
=
| \varphi |^2 (N_+^2)_{ij}
\eea
and, using (\ref{contract}), we compute
\bea
({\rm II})\cdot\varphi_{iab}\o^{ka}\o^{pb}
&=&
{| \varphi|^2 \over 4}g^{\ell m}
(\o_{\lambda i}g_{\mu a}+\o_{\mu a}g_{\lambda i}
-
\o_{\lambda a}g_{\mu i}-\o_{\mu i} g_{\lambda a})N_{mk}{}^\mu N_{\ell j}{}^\lambda\o^{ka}
\nonumber\\
&=&
{| \varphi|^2 \over 4}g^{\ell m}
(\o_{\lambda i}J^k{}_\mu-\delta^k{}_\mu g_{\lambda i}+\delta^k{}_\lambda g_{\mu i}
-
\o_{\mu i}J^k{}_\lambda)N_{mk}{}^\mu N_{\ell j}{}^\lambda
\nonumber\\
&=&
{| \varphi|^2 \over 4}
(\o_{\lambda i}N_{mk}{}^{Jk}N^m{}_j{}^\lambda
-
N_{mk}{}^kN^m{}_{j}{}^\lambda
+
N_{mki}N^m{}_j{}^k
-
\o_{\mu i}N_{mk}{}^\mu N^m{}_{j}{}^{Jk}).
\nonumber
\eea
Now $N_{mk}{}^k=0$, and by the Nijenhuis tensor identities, 
\bea
N_{mk}{}^{Jk}=N_{m,Jk}{}^k=N_{Jm,k}{}^k=0.
\eea
Furthermore, we have by definition $
N_{mki}N^m{}_j{}^k=-(N_+^2)_{ij}$, while
\bea
-
\o_{\mu i}N_{mk}{}^\mu N^m{}_{j}{}^{Jk}
&=&
g_{\mu \nu}J^\nu{}_iN_{mk}{}^\mu N^m{}_{j}{}^{Jk}
=
N_{mk\,Ji}N^m{}_j{}^{Jk}
=
N_{m,Jk,i}N^m{}_j{}^{Jk}
\nonumber\\
&=&
-N_{mki}N^m{}_j{}^k=(N_+^2)_{ij}
\eea
and hence
\bea
({\rm II})\cdot\varphi_{iab}\o^{ka}\o^{pb}=0.
\eea
Since (${\rm III}$) can be obtained from (${\rm II}$) by the simultaneous interchange $a\leftrightarrow b$ and $k\leftrightarrow p$, we also have
\bea
({\rm III})\cdot\varphi_{iab}\o^{ka}\o^{pb}=0.
\eea

\medskip
We consider next the expression
\bea
g^{\ell m}E_{m;j\lambda p}N_{\ell k}{}^\lambda \varphi_{iab}\o^{ka}\o^{pb}
&=&
g^{\ell m}(\varphi_{\mu\lambda p}N_{mj}{}^\mu
+
\varphi_{j\mu p}N_{m\lambda}{}^\mu
+
\varphi_{j\lambda\mu}N_{mp}{}^\mu)
N_{\ell k}{}^\lambda \varphi_{iab}\o^{ka}\o^{pb}
\nonumber\\
&=&
({\rm IV}+{\rm V}+{\rm VI})
\varphi_{iab}\o^{ka}\o^{pb}.
\eea
The contributions of the term (IV) worked out to be $0$,
\bea
({\rm IV})
\cdot\varphi_{iab}\o^{ka}\o^{pb}
&=&
{| \varphi|^2 \over 4}
(\o_{\mu i}g_{\lambda a}+\o_{\lambda a}g_{\mu i}-\o_{\mu a}g_{\lambda i}-\o_{\lambda i}g_{\mu a})
N_{mj}{}^\mu N^m{}_k{}^\lambda\o^{ka}
\nonumber\\
&=&
{| \varphi|^2 \over 4}
(\o_{\mu i}J^k{}_\lambda-\delta^k{}_\lambda g_{\mu i}+\delta^k{}_\mu g_{\lambda i}
-\o_{\lambda i}J^k{}_\mu)N_{mj}{}^\mu N^m{}_k{}^\lambda.
\eea
The first two terms on the right hand side vanish individually, since
\bea
&&
\o_{\mu i}J^k{}_\lambda N_{mj}{}^\mu N^m{}_k{}^\lambda
=
\o_{\mu i}N_{mj}{}^\mu N^m{}_{J\lambda}{}^\lambda=0
\nonumber\\
&&
\delta^k{}_\lambda \,g_{\mu i}
N_{mj}{}^\mu N^m{}_k{}^\lambda
=
N_{mj}{}^\mu N^m{}_k{}^k=0.
\eea
Of the remaining two terms, we have obviously
\bea
\delta^k{}_\mu g_{\lambda i}
N_{mj}{}^\mu N^m{}_k{}^\lambda
=
N_{mj}{}^kN^m{}_{ki}
=
-N_m{}^k{}_jN^m{}_{ki}=-(N_+^2)_{ij},
\eea
while
\bea
-\o_{\lambda i}N_{mj}{}^{Jk}N^m{}_k{}^\lambda
&=&
g_{\lambda \nu}J^\nu{}_iN_{mj}{}^{Jk}N^m{}_k{}^\lambda
=
J^\nu{}_iN_{mj}{}^{Jk}N^m{}_{k\nu}
=
N_{mj}{}^{Jk}N^m{}_{k,Ji}
\nonumber\\
&=&
N_{mj}{}^{Jk}N^m{}_{Jk,i}
=
-N_{mj}{}^kN^m{}_{ki}
=
(N_+^2)_{ij}
\eea
so they cancel each other out and we obtain, as claimed,
\bea
({\rm IV})
\cdot\varphi_{iab}\o^{ka}\o^{pb}=0.
\eea
The next group of terms is given by
\bea
({\rm V})\cdot
\varphi_{iab}\o^{ka}\o^{pb}
&=&
\varphi_{j\mu p}N_{m\lambda}{}^\mu N^m{}_k{}^\lambda
\varphi_{iab}\o^{ka}\o^{pb}=
-\varphi_{j\mu p}g^{\mu\nu}(N_+^2)_{\nu k}\varphi_{iab}\o^{ka}\o^{pb}
\nonumber\\
&=&
-
{| \varphi|^2 \over 4}(\o_{ji}g_{\mu a}+\o_{\mu a}g_{ji}-\o_{ja}g_{\mu i}-\o_{\mu i}g_{ja})
\o^{ka}(N_+^2)_{\nu k}g^{\mu\nu}.
\eea
The first term on the right produces $0$, since it can be computed as $\o_{ji} \o^{k\nu} (N_+^2)_{\nu k}$. This term vanishes due to the anti-symmetrization of $k$ and $\nu$. 
We are left with
\bea
({\rm V})\cdot
\varphi_{iab}\o^{ka}\o^{pb}
&=&
-{| \varphi|^2 \over 4}
(-\delta^k{}_\mu g_{ji}g^{\mu \nu}
+
\delta^k{}_j g_{\mu i}g^{\mu \nu}
-
\o_{\mu i}J^k{}_j g^{\mu \nu})(N_+^2)_{\nu k}
\nonumber\\
&=&
-{| \varphi|^2 \over 4}(-|N|^2 g_{ij}+(N_+^2)_{ij}
+(N_+^2)_{Ji ,Jj})
\eea
Since we have
\bea
(N_+^2)_{Ji ,Jj}=N^{mk}{}_{Ji}N_{mk,Jj}
=
- N^{m,Jk}{}_iN_{m,Jk,j}
=
N^{mk}{}_iN_{mkj}=(N_+^2)_{ij}
\eea
we are left with
\bea
({\rm V})\cdot
\varphi_{iab}\o^{ka}\o^{pb}
=
{| \varphi|^2 \over 4}|N|^2 g_{ij} - {| \varphi |^2 \over 2} (N_+^2)_{ij}.
\eea
Finally, we observe that
\bea
({\rm VI})\cdot
\varphi_{iab}\o^{ka}\o^{pb}=
g^{\ell m}
\varphi_{j\lambda\mu}N_{mp}{}^\mu
N_{\ell k}{}^\lambda \varphi_{iab}\o^{ka}\o^{pb}
=0.
\eea
We can readily see this in a complex frame. Since $\varphi\in \Lambda^{3,0}\oplus\Lambda^{0,3}$, the only components of $\varphi_{j\lambda\mu}$ which are not $0$ must have both barred or both unbarred indices. But the contraction with $g^{\ell m}$ implies that the indices $\ell $ and $m$ must be mixed. But then for $N_{mp}{}^\mu
N_{\ell k}{}^\lambda$ not to be $0$, the indices $\lambda$ and $\mu$ must be mixed too, contradicting the requirement that they must be both barred or both unbarred. This establishes our claim.

\medskip
We still have one more contribution from the second row of (\ref{NaE}), given by
\bea
g^{\ell m}E_{m;jk\lambda}N_{\ell p}{}^\lambda\varphi_{iab}\o^{ka}\o^{pb}
\eea
but which can be recognized as coinciding with the term that we just computed
\bea
g^{\ell m}E_{m;j\lambda p}N_{\ell k}{}^\lambda\varphi_{iab}\o^{ka}\o^{pb}={ | \varphi|^2 \over 4}|N|^2 g_{ij}- {| \varphi |^2 \over 2} (N_+^2)_{ij}
\eea
upon the renaming of indices $a\leftrightarrow b$, $p\leftrightarrow k$. 

\medskip
It is convenient to summarize the formula which we have obtained as a lemma:

\begin{lemma}
We have
\bea
g^{\ell m}(E_{m;\lambda kp}N_{\ell j}{}^\lambda
+
E_{m;j\lambda p} N_{\ell k}{}^\lambda
+
E_{m;jk\lambda}N_{\ell p}{}^\lambda)
\varphi_{iab}\o^{ka}\o^{pb}
=
{| \varphi|^2 \over 2}|N|^2 g_{ij}.
\eea
\end{lemma}

\subsubsection{Contributions of the term $E_{m;\ell jkp}$}

The term $E_{m;\ell jkp}$ involves ${\frak D}_\mu\varphi_{jkp}$, ${\frak D}_\ell\varphi_{\mu kp}$, ${\frak D}_\ell\varphi_{j\mu p}$, and ${\frak D}_\ell\varphi_{jk\mu}$. We use (\ref{frakD-varphi}) to evaluate the contribution of each term in turn,
\bea
{\frak D}_\mu\varphi_{jkp}
N_{m\ell}{}^\mu \varphi_{iab}\o^{ka}\o^{pb}
&=&
-{1\over 2}(\alpha_\mu\varphi_{jkp}+\alpha_{J\mu}\hat\varphi_{jkp})N_{m\ell}{}^\mu \varphi_{iab}\o^{ka}\o^{pb}
\nonumber\\
&=&
{1\over 2}| \varphi|^2 (\alpha_\mu g_{ji}-\alpha_{J\mu}\o_{ji})N_{m\ell}{}^\mu
\eea
Upon symmetrization in $i$ and $j$, and contracting with $g^{\ell m}$, we obtain
\bea
g^{\ell m}
{\frak D}_\mu\varphi_{jkp}
N_{m\ell}{}^\mu \varphi_{iab}\o^{ka}\o^{pb}+(i\leftrightarrow j)
=
| \varphi|^2 \, g_{ij} \alpha_\mu g^{m\ell}N_{m\ell}{}^\mu=0
\eea
where we have used the fact that $N$ is of type $(0, 2)$ to write
\bea
g^{m\ell}N_{m\ell}{}^\mu=g^{Jm, J\ell} N_{m\ell}{}^\mu = g^{m \ell} N_{Jm, J\ell}{}^\mu = - g^{m\ell} N_{m\ell}{}^\mu\nonumber
\eea
and therefore
\bea
g^{m\ell}N_{m\ell}{}^\mu=0.
\eea

Next, we consider the term
\bea
{\frak D}_\ell\varphi_{\mu kp}N_{mj}{}^\mu\varphi_{iab}\o^{ka}\o^{pb}
&=&
-{1\over 2}(\alpha_\ell\varphi_{\mu kp}+\alpha_{J\ell}\hat\varphi_{\mu kp})
N_{mj}{}^\mu\varphi_{iab}\o^{ka}\o^{pb}\\
\nonumber
&=&
{1\over 2}| \varphi|^2 (\alpha_\ell g_{\mu i}-\alpha_{J\ell}\o_{\mu i})N_{mj}{}^\mu
=
{1\over 2}| \varphi|^2 (\alpha_\ell N_{mji}+\alpha_{J\ell}N_{mj,Ji}).
\eea
The first term on the right symmetrizes to $0$. So does the second, using the fact that $N$ is a type $(0,2)$-tensor, so that
$N_{mj,Ji}=N_{Jm,j,i}$ which is antisymmetric in the last two indices.

We consider now
\bea
{\frak D}_\ell\varphi_{j\mu p}
N_{mk}{}^\mu\varphi_{iab}\o^{ka}\o^{pb}
&=&
-
{1\over 2}(\alpha_\ell\varphi_{j\mu p}+\alpha_{J\ell}\hat\varphi_{j\mu p})
N_{mk}{}^\mu \varphi_{iab}\o^{ka}\o^{pb}.
\eea
We work out separately the contributions  of the two terms $\varphi_{j\mu p}$ and $\hat\varphi_{j\mu p}$ on the right hand side. First, we have
\bea
\alpha_\ell\varphi_{j\mu p}
N_{mk}{}^\mu \varphi_{iab}\o^{ka}\o^{pb}
&=&
{| \varphi|^2 \over 4}\alpha_\ell
(\o_{ji}g_{\mu a}+\o_{\mu a}g_{ji}-\o_{ja}g_{\mu i}-\o_{\mu i}g_{ja})
N_{mk}{}^\mu \o^{ka}
\nonumber\\
&=&
{| \varphi|^2 \over 4}\alpha_\ell (\o_{ji}J^k{}_\mu
-\delta^k{}_\mu g_{ji}+\delta^k{}_jg_{\mu i}-\o_{\mu i}J^k{}_j)N_{mk}{}^\mu.
\eea
We claim that, upon symmetrization in $i$ and $j$, the net result is $0$. This is obviously true of the term $\o_{ji}$, while $N_{mk}{}^k=0$ and $N_{mji}$ also symmetrizes to $0$. The fourth term can be rewritten as
\bea
\o_{\mu i}J^k{}_jN_{mk}{}^\mu
=
-g_{\mu \nu}J^\nu{}_i N_{m,Jj}{}^\mu
=
-J^\nu{}_iN_{m,Jj,\nu}=N_{m,Jj,Ji}
\eea
which symmetrizes to $0$. We come to the contribution of the term involving $\hat\varphi$,
\bea
\hat\varphi_{j\mu p}
\varphi_{iab}\o^{ka}\o^{pb}
&=&
-\varphi_{Jj,\mu p}\varphi_{iab}\o^{ka}\o^{pb}
\nonumber\\
&=&-{| \varphi|^2 \over 4}
(\o_{Jj,i}g_{\mu a}+\o_{\mu a}g_{Jj,i}-\o_{Jj,a}g_{\mu i}-\o_{\mu i}g_{Jj,a})\o^{ka}
\nonumber\\
&=&
-{| \varphi|^2 \over 4}(-g_{ij}g_{\mu a}+\o_{\mu a}\o_{ji}
+
g_{aj}g_{\mu i}-\o_{\mu i}\o_{ja})\o^{ka}
\nonumber\\
&=&
-{| \varphi|^2 \over 4}(-g_{ij}J^k{}_\mu
-
\delta^k{}_\mu \o_{ji}+J^k{}_j g_{\mu i}+\delta^k{}_j\o_{\mu i}).
\eea
Dropping the term $\o_{ji}$ since it symmetrizes to $0$, we arrive at
\bea
\hat\varphi_{j\mu p}
\varphi_{iab}\o^{ka}\o^{pb}N_{mk}{}^\mu+(i\leftrightarrow j)
&=&
-{| \varphi|^2 \over 4}
(-g_{ij}N_{m,J\mu}{}^\mu
+J^k{}_jN_{mki}
+N_{mj}{}^\mu\o_{\mu i})
+(i\leftrightarrow j)
\nonumber\\
&=&
-{| \varphi|^2 \over 4}(N_{m,Jj,i}-N_{mj,Ji})+
(i\leftrightarrow j)
\nonumber\\
&=&
-{| \varphi|^2 \over 4}(N_{Jm,j,i}-N_{Jm,j,i})+(i\leftrightarrow j)=0.
\eea
The last term ${\frak D}_\ell\varphi_{jk\mu}$ makes an identical contribution as ${\frak D}_{\ell}\varphi_{ j\mu p}$, upon renaming the summation indices $a\leftrightarrow b$, $k\leftrightarrow p$. Thus its contribution is also $0$. In summary, we have established

\begin{lemma}
We have
\bea
g^{m\ell}
E_{m;\ell jkp}\varphi_{iab}\o^{ka}\o^{pb}+(i\leftarrow j)
=0.
\eea
\end{lemma}

\subsubsection{Completion of the calculations for $\na^\ell E_{\ell;jkp}$}

The terms from $\na^\ell E_{\ell;jkp}$ in (\ref{NaE}) whose contributions we have not worked out as yet are the following
\bea
&&
(g^{\ell m}{\frak D}_m\varphi_{\lambda kp})N_{\ell j}{}^\lambda
+
(g^{\ell m}{\frak D}_m\varphi_{j\lambda p})N_{\ell k}{}^\lambda
+
(g^{\ell m}{\frak D}_m\varphi_{jk\lambda})N_{\ell p}{}^\lambda
\nonumber\\
&&
+
\varphi_{\lambda kp}g^{\ell m}\na_mN_{\ell j}{}^\lambda
+
\varphi_{j\lambda p}g^{\ell m}\na_mN_{\ell k}{}^\lambda
+
\varphi_{jk\lambda}g^{\ell m}\na_mN_{\ell p}{}^\lambda.
\nonumber\\
&&
=
\varphi_{\lambda kp}(-{1\over 2}\alpha^\ell N_{\ell j}{}^\lambda+\na^\ell N_{\ell j}{}^\lambda)
-
{1\over 2}\hat\varphi_{\lambda kp}\,\alpha_{Jm}g^{\ell m}N_{\ell j}{}^\lambda
\nonumber\\
&&
\ +
\varphi_{j\lambda p}(-{1\over 2}\alpha^\ell N_{\ell k}{}^\lambda+\na^\ell N_{\ell k}{}^\lambda)
-
{1\over 2}\hat\varphi_{j\lambda p}\,\alpha_{Jm}g^{\ell m}N_{\ell k}{}^\lambda
\nonumber\\
&&
\ +
\varphi_{jk\lambda}(-{1\over 2}\alpha^\ell N_{\ell p}{}^\lambda+\na^\ell N_{\ell p}{}^\lambda)
-{1\over 2}\hat\varphi_{jk\lambda}\,\alpha_{Jm}g^{\ell m}N_{\ell p}{}^\lambda
\nonumber\\
&&
=\
{\rm VII}+\hat{\rm VII}
+
{\rm VIII}+\hat{\rm VIII}+{\rm IX}+\hat{\rm IX}.
\eea
Again, we evaluate each contribution in turn. We have
\bea
({\rm VII})\cdot
\varphi_{iab}\o^{ka}\o^{pb}
=
-| \varphi|^2 g_{\lambda i}(-{1\over 2}\alpha^\ell N_{\ell j}{}^\lambda+\na^\ell N_{\ell j}{}^\lambda)
=
| \varphi|^2 ({1\over 2}\alpha^\ell N_{\ell ji}-\na^\ell N_{\ell ji}) =0
\eea
upon symmetrization in $i\leftrightarrow j$. Next,
\bea
(\hat{\rm VII})\cdot
\varphi_{iab}\o^{ka}\o^{pb}
&=&
-{1\over 2}\alpha_{Jm}\varphi_{J\lambda,k,p}N^m{}_j{}^\lambda \varphi_{iab}\o^{ka}\o^{pb}
\nonumber\\
&=&
{1\over 2}| \varphi|^2 \alpha_{Jm}\o_{\lambda i}N^m{}_j{}^\lambda
=
-{1\over 2}| \varphi|^2 \alpha_{Jm}N^m{}_{j,Ji}\nonumber\\
&=&-{1\over 2}| \varphi|^2 \alpha_{Jm}g^{m\ell}N_{\ell ,j ,Ji}
=-{1\over 2}| \varphi|^2 \alpha_{Jm}g^{m\ell}N_{J\ell,j,i}
\eea
which produces $0$ upon symmetrization in $j$ and $i$. Next,
\bea
({\rm VIII})\cdot\varphi_{iab}\o^{ka}\o^{pb}
&=&
\varphi_{j\lambda p}(-{1\over 2}\alpha^\ell N_{\ell k}{}^\lambda+\na^\ell N_{\ell k}{}^\lambda)\varphi_{iab}
\o^{ka}\o^{pb}
\nonumber\\
&=&
{| \varphi|^2 \over 4}(\o_{ji}g_{\lambda a}+\o_{\lambda a}g_{ji}-\o_{ja}g_{\lambda i}-\o_{\lambda i}g_{ja})\o^{ka}(-{1\over 2}\alpha^\ell N_{\ell k}{}^\lambda+\na^\ell N_{\ell k}{}^\lambda)
\nonumber\\
&=&
{| \varphi|^2 \over 4}(\o_{ji}J^k{}_\lambda
-\delta^k{}_\lambda g_{ji}+\delta^k{}_j g_{\lambda i}-\o_{\lambda i}J^k{}_j)
(-{1\over 2}\alpha^\ell N_{\ell k}{}^\lambda+\na^\ell N_{\ell k}{}^\lambda)
\nonumber\\
&=&
{| \varphi|^2 \over 4}
\big\{g_{ji}({1\over 2}\alpha^\ell N_{\ell\lambda}{}^\lambda-\na^\ell N_{\ell\lambda}{}^\lambda)
-{1\over 2}\alpha^\ell N_{\ell ji}+\na^\ell N_{\ell ji}
-{1\over 2}\alpha^\ell N_{\ell, Jj,Ji}\nonumber\\
&& - \o_{\lambda i} J^k{}_j \nabla^\ell N_{\ell k}{}^\lambda\big\}
\eea
Note that, the first two terms are zero because $N_{\ell \lambda}{}^\lambda =0$; the next three terms also adds up to $0$ upon symmetrization in $i$ and $j$. Indeed, the last term is also zero upon symmetrization in $i$ and $j$ because it is anti-symmetric about $i$ and $j$ as
\bea
\o_{\lambda i} J^k{}_j \nabla^\ell N_{\ell k}{}^\lambda=\o_{\lambda i} g^{k p} J_{pj} \nabla^\ell N_{\ell k}{}^\lambda = \o_{\lambda i }  \o_{pj} \nabla^\ell (N_{\ell k}{}^\lambda g^{kp})
=
\o_{\lambda i} \o_{pj} \nabla^\ell N_\ell{}^{p\lambda}= - \o_{\lambda j} \o_{p i} \nabla^\ell N_\ell{}^{p\lambda} \nonumber
\eea
The last identity is seen by switching indices $p\leftrightarrow \lambda$ and using the antisymmetry of $N$.

The next term to be considered is
\bea
(\hat{\rm VIII})\cdot \varphi_{iab}\o^{ka}\o^{pb}
&=&
-{1\over 2}\alpha_{Jm}\hat\varphi_{j\lambda p}
g^{\ell m}N_{\ell k}{}^\lambda
\varphi_{iab}\o^{ka}\o^{pb}
=
{1\over 2}\alpha_{Jm}\varphi_{Jj,\lambda,p}
g^{\ell m}N_{\ell k}{}^\lambda
\varphi_{iab}\o^{ka}\o^{pb}
\nonumber\\
&=&
{| \varphi|^2 \over 8}(\o_{Jj,i}g_{\lambda a}+\o_{\lambda a}g_{Jj,i}
-
\o_{Jj,a}g_{\lambda i}-\o_{\lambda i}g_{Jj,a})\o^{ka}\alpha_{Jm}N^m{}_k{}^\lambda
\nonumber\\
&=&
{| \varphi|^2 \over 8}(-g_{ij}J^k{}_\lambda-\delta^k{}_\lambda \o_{ji}
+J^k{}_j g_{\lambda i}+\o_{\lambda i}\delta^k{}_j)
\alpha_{Jm}N^m{}_k{}^\lambda
\nonumber\\
&=&
{| \varphi|^2 \over 8}
(-g_{ij}\alpha_{Jm}N^m{}_{J\lambda}{}^\lambda-\alpha_{Jm}N^m{}_\lambda{}^\lambda
+\alpha_{Jm}N^m{}_{Jj,i}-\alpha_{Jm}N^m{}_{j,Ji})
\nonumber
\eea
Using the fact that $N$ is a tensor of type $(0,2)$, we readily see that each of these terms reduces to $0$.

\medskip
In summary, the contribution of the Laplacian term is given by

\begin{lemma} \label{laplacian-contribution}
We have
\bea 
(g^{\ell m}\na_m\na_\ell\varphi_{jkp})\varphi_{iab}\o^{ka}\o^{pb}
+
(i\leftrightarrow j)
=
| \varphi|^2 \big\{ \nabla_\mu \alpha^\mu  + |N|^2\big\}g_{ij}.
\nonumber
\eea
\end{lemma}

\medskip

\subsubsection{Contributions of the curvature terms}

Turning next to the curvature contributions, we write
\bea
g^{\ell m}
[\na_m,\na_j]\varphi_{kp\ell}
&=&
-g^{\ell m}(R_{mj}{}^\lambda{}_k\varphi_{\lambda p\ell}
+
R_{mj}{}^\lambda{}_p\varphi_{k\lambda\ell}
+
R_{mj}{}^\lambda{}_\ell\varphi_{kp\lambda})
\nonumber\\
&=&
-R^\ell{}_j{}^{\lambda}{}_k\varphi_{\lambda p\ell}-R^\ell{}_j{}^{\lambda}{}_p\varphi_{k\lambda\ell}+R_j{}^\lambda\varphi_{kp\lambda}
\nonumber\\
&=&
-R^\ell{}_j{}^{\lambda}{}_k\varphi_{\lambda p\ell}+R^\ell{}_j{}^{\lambda}{}_p\varphi_{\lambda k\ell}+R_j{}^\lambda\varphi_{kp\lambda}
\eea
We consider for the moment only the contribution of the last term.
\bea
R_j{}^\lambda\,\varphi_{kp\lambda}
\,\varphi_{iab}\o^{ka}\o^{pb}
=
- | \varphi |^2 R_j{}^\lambda \,g_{\lambda i}=-| \varphi |^2 R_{ji}.
\eea
The next curvature contribution is similar
\bea
g^{\ell m}[\na_m,\na_p]\varphi_{jk\ell}
&=&
-g^{\ell m}(R_{mp}{}^\lambda{}_j\varphi_{\lambda k\ell}
+
R_{mp}{}^\lambda{}_k\varphi_{j\lambda\ell}
+
R_{mp}{}^\lambda{}_\ell\varphi_{jk\lambda})
\nonumber\\
&=&
-R^\ell{}_p{}^{\lambda}{}_j\varphi_{\lambda k\ell}+R^\ell{}_p{}^{\lambda}{}_k\varphi_{\lambda j\ell}
+R_p{}^\lambda\varphi_{jk\lambda}
\eea
and the corresponding last term gives
\bea
R_p{}^\lambda\varphi_{jk\lambda}\,\varphi_{iab}\o^{ka}\o^{pb}
&=&
{| \varphi |^2 \over 4}R_p{}^\lambda(\o_{ji}g_{\lambda b}+\o_{\lambda b}g_{ji}
-\o_{jb}g_{\lambda i}-\o_{\lambda i}g_{jb})\o^{pb}
\nonumber\\
&=&
{| \varphi |^2 \over 4}R_p{}^\lambda
(\o_{ji}J^p{}_\lambda-\delta^p{}_\lambda g_{ji}+\delta^p{}_j g_{\lambda i}-\o_{\lambda i}J^p{}_j)
\nonumber\\
&=&
{| \varphi |^2 \over 4}(-R\,g_{ji}+R_{ji}+R_{Jj,Ji})
\eea
where we have dropped the term proportional to $\o_{ji}$ since it symmetrizes to $0$. The remaining terms gives an identical contribution. Indeed,
\bea
g^{\ell m}[\na_m,\na_k]\varphi_{pj \ell}
&=&
-g^{\ell m}(R_{mk}{}^\lambda{}_p \varphi_{\lambda j\ell}
+
R_{mk}{}^\lambda{}_j\varphi_{p\lambda\ell}
+
R_{mk}{}^\lambda{}_\ell\varphi_{pj \lambda})
\nonumber\\
&=&
-R^\ell{}_k{}^{\lambda}{}_p \varphi_{\lambda j \ell}+R^\ell{}_k{}^{\lambda}{}_j\varphi_{\lambda p\ell}+R_k{}^\lambda\,\varphi_{pj \lambda}.
\eea
Considering for the moment only the contribution of the last term, we can write
\bea
R_k{}^\lambda\,\varphi_{pj \lambda}\varphi_{iab}\o^{ka}\o^{pb}
&=&
R_k{}^\lambda (\o^{pb} \varphi_{pj\lambda}\varphi_{bia}) \o^{ka}
\nonumber\\
&=&
{| \varphi|^2 \over 4} R_k{}^\lambda(\o_{ji}g_{\lambda a}+\o_{\lambda a}g_{ji}-\o_{ja}g_{\lambda i}-\o_{\lambda i}g_{ja})\o^{ka}
\nonumber\\
&=&
{| \varphi|^2 \over 4}R_k{}^\lambda(\o_{ji}J^k{}_\lambda-\delta^k{}_\lambda g_{ji}+\delta^k{}_jg_{\lambda i}
-\o_{\lambda i}J^k{}_j)
\nonumber\\
&=&
{| \varphi|^2 \over 4}(-R\,g_{ji}+R_{ji}+R_{Jj,Ji})
\eea
where we have dropped the antisymmetric term $
\o_{ji}$ just as before. Assembling all the terms, we have proved the following lemma

\begin{lemma}
We have the following formula
\bea
&&
g^{\ell m}([\na_m,\na_j]\varphi_{kp\ell}
+
[\na_m,\na_k]\varphi_{pj\ell}+[\na_m,\na_p]\varphi_{jk\ell})
\varphi_{iab}\o^{ka}\o^{pb}+(i\leftrightarrow j)
\nonumber\\
&&
=-2 | \varphi |^2 \,R_{ji} - | \varphi |^2 R\, g_{ij} + | \varphi |^2 R_{ij} + | \varphi |^2 R_{Jj,Ji} +F
\eea
where the term $F$ is given by
\bea
F&=&\big\{(R^\ell{}_j{}^{\lambda}{}_p-R^\ell{}_p{}^{\lambda}{}_j)\varphi_{\lambda k\ell}
+(-R^\ell{}_j{}^{\lambda}{}_k+R^\ell{}_k{}^{\lambda}{}_j)\varphi_{\lambda p\ell}
+(R^\ell{}_p{}^{\lambda}{}_k-R^\ell{}_k{}^{\lambda}{}_p)\varphi_{\lambda j\ell}\big\}
\varphi_{iab}\o^{ka}\o^{pb}
\nonumber\\
&&
+(i\leftrightarrow j)
\eea
\end{lemma}

\subsubsection{Evaluation of the term $F$}

We begin with
\bea
(R^\ell{}_j{}^{\lambda}{}_p-R^\ell{}_p{}^{\lambda}{}_j)\varphi_{\lambda k\ell}
\varphi_{iab}\o^{ka}\o^{pb}
&=&
{| \varphi|^2 \over 4}
(R^\ell{}_j{}^{\lambda}{}_p-R^\ell{}_p{}^{\lambda}{}_j)
(\o_{\lambda i}g_{\ell b}+\o_{\ell b}g_{\lambda i}-\o_{\lambda b}g_{\ell i}
-\o_{\ell i}g_{\lambda b})\o^{pb}
\nonumber\\
&=&
{| \varphi|^2 \over 4}
(R^\ell{}_j{}^{\lambda}{}_p-R^\ell{}_p{}^{\lambda}{}_j)
(\o_{\lambda i}J^p{}_\ell-\delta^p{}_\ell g_{\lambda i}+\delta^p{}_\lambda g_{\ell i}-\o_{\ell i}J^p{}_\lambda)
\nonumber\\
&=&
{| \varphi|^2 \over 4}(-R^{Jp}{}_{j,Ji,p}+R^\ell{}_{J\ell}{}_{Ji,j})
-{| \varphi|^2 \over 4}(-R_{ji}+R^\ell{}_{\ell ij})
\nonumber\\
&&
+
{| \varphi|^2 \over 4}(R_{ij}{}^\lambda{}_\lambda-R_{i \lambda}{}^\lambda{}_j)
+{| \varphi|^2 \over 4}(R_{Ji,j}{}^\lambda{}_{J\lambda}-R_{Ji,J\lambda}{}^\lambda{}_j)
\eea
This reduces to
\bea
(R^\ell{}_j{}^{\lambda}{}_p-R^\ell{}_p{}^{\lambda}{}_j)\varphi_{\lambda k\ell}
\varphi_{iab}\o^{ka}\o^{pb}
=
{| \varphi|^2 \over 4}(-R{}^{Jp}{}_{j,Ji,p}+R^\ell{}_{J
  \ell, Ji,j}
+R_{Ji,j}{}^\lambda{}_{J\lambda}-R_{Ji,J\lambda}{}^\lambda{}_j)
+{| \varphi|^2 \over 2}R_{ij}.
\nonumber
\eea
Using the symmetries of the Riemann curvature tensor, we simplify this to
\bea
\label{F1}
(R^\ell{}_j{}^{\lambda}{}_p-R^\ell{}_p{}^{\lambda}{}_j)\varphi_{\lambda k\ell}
\varphi_{iab}\o^{ka}\o^{pb}
=
{| \varphi|^2 \over 2}(-R_{Ji,J \lambda}{}^\lambda{}_j +R_{Ji,j}{}^\lambda{}_{J\lambda})
+{| \varphi|^2 \over 2}R_{ij}.
\nonumber
\eea
We work out the next term, which after relabeling is
\bea
\label{F2}
&&(-R^\ell{}_j{}^{\lambda}{}_k+R^\ell{}_k{}^{\lambda}{}_j)\varphi_{\lambda p\ell}
\varphi_{iab}\o^{ka}\o^{pb}\nonumber\\
&=&
(R^\ell{}_j{}^{\lambda}{}_k-R^\ell{}_k{}^{\lambda}{}_j)\varphi_{\lambda \ell p}
\varphi_{iab}\o^{ka}\o^{pb}
=
(R^\ell{}_j{}^{\lambda}{}_p-R^\ell{}_p{}^{\lambda}{}_j)\varphi_{\lambda \ell k}
\varphi_{iba}\o^{pb}\o^{ka},
\eea
and is therefore identical to the previous term,
\be
(-R^\ell{}_j{}^{\lambda}{}_k+R^\ell{}_k{}^{\lambda}{}_j)\varphi_{\lambda p\ell}
\varphi_{iab}\o^{ka}\o^{pb} = {| \varphi|^2 \over 2}(-R_{Ji,J \lambda}{}^\lambda{}_j +R_{Ji,j}{}^\lambda{}_{J\lambda})
+{| \varphi|^2 \over 2}R_{ij}.
\ee
We work out the final term. We start with
\be
(R^\ell{}_p{}^{\lambda}{}_k - R^\ell{}_k{}^{\lambda}{}_p)\varphi_{\lambda j \ell} \varphi_{iab}\o^{ka}\o^{pb} = -R_{pk}{}^{\lambda \ell} \varphi_{\lambda j \ell} \varphi_{iab}\o^{ka}\o^{pb}
\ee
by the Bianchi identity $R^\ell{}_p{}^{\lambda}{}_k + R_{pk}{}^{\lambda \ell} + R_k{}^\ell{}^{\lambda}{}_p =0$.
Applying the identity (\ref{riemann-J-identity}) gives
\bea
-R_{pk}{}^{\lambda \ell} \varphi_{\lambda j \ell} \varphi_{iab}\o^{ka}\o^{pb} &=& (-R_{pk}{}^{J\lambda, J\ell} + B_{kp}{}^{\lambda \ell} ) \varphi_{\lambda j \ell} \varphi_{iab}\o^{ka}\o^{pb} \nonumber\\
&=& - R_{pk}{}^{\lambda \ell} \varphi_{J \lambda, j, J\ell} \varphi_{iab}\o^{ka}\o^{pb} + B_{kp}{}^{\lambda \ell} \varphi_{\lambda j \ell} \varphi_{iab}\o^{ka}\o^{pb} \nonumber\\
&=&  R_{pk}{}^{\lambda \ell} \varphi_{ \lambda, j, \ell} \varphi_{iab}\o^{ka}\o^{pb} + B_{kp}{}^{\lambda \ell} \varphi_{\lambda j \ell} \varphi_{iab}\o^{ka}\o^{pb}
\eea
Therefore
\be
-R_{pk}{}^{\lambda \ell} \varphi_{\lambda j \ell} \varphi_{iab}\o^{ka}\o^{pb} =  {1 \over 2} B_{kp}{}^{\lambda \ell}  \varphi_{\lambda j \ell} \varphi_{iab}\o^{ka}\o^{pb}
\ee
and hence
\be
(R^\ell{}_p{}^{\lambda}{}_k - R^\ell{}_k{}^{\lambda}{}_p)\varphi_{\lambda j \ell} \varphi_{iab}\o^{ka}\o^{pb} = {1 \over 2} B_{kp}{}^{\lambda \ell} \varphi_{\lambda j \ell} \varphi_{iab}\o^{ka}\o^{pb}
\ee
By definition of $B$,
\be
B_{kp}{}^{\lambda \ell} = - 2 {\frak D}_k N_{p}{}^{\lambda \ell} + 2 {\frak D}_p N_k{}^{\lambda \ell} - 2 N^\alpha{}_{kp} N_{\alpha}{}^{\lambda \ell} .
\ee
Therefore
\bea \label{F-terms}
(R^\ell{}_p{}^{\lambda}{}_k - R^\ell{}_k{}^{\lambda}{}_p)\varphi_{\lambda j \ell} \varphi_{iab}\o^{ka}\o^{pb} 
&=& - {\frak D}_k N_{p}{}^{\lambda \ell} \varphi_{\lambda j \ell} \varphi_{iab}\o^{ka}\o^{pb} + {\frak D}_p N_k{}^{\lambda \ell} \varphi_{\lambda j \ell} \varphi_{iab}\o^{ka}\o^{pb} \nonumber\\
&&-  N^\alpha{}_{kp} N_{\alpha}{}^{\lambda \ell}\varphi_{\lambda j \ell} \varphi_{iab}\o^{ka}\o^{pb} \\
&=&-2 {\frak D}_k N_{p}{}^{\lambda \ell} \varphi_{\lambda j \ell} \varphi_{iab}\o^{ka}\o^{pb} -  N^\alpha{}_{kp} N_{\alpha}{}^{\lambda \ell}\varphi_{\lambda j \ell} \varphi_{iab}\o^{ka}\o^{pb}.
\nonumber
\eea
We start with the last term. By the Bianchi identity $N_{ijk} + N_{kij} + N_{jki} = 0$,
\be
N^\alpha{}_{kp} N_{\alpha}{}^{\lambda \ell}\varphi_{\lambda j \ell} \varphi_{iab}\o^{ka}\o^{pb} = - N^\alpha{}_{kp} [ N^\ell{}_{\alpha}{}^{\lambda}\varphi_{\lambda j \ell} + N^{\lambda \ell}{}_\alpha \varphi_{\lambda j \ell}] \varphi_{iab}\o^{ka}\o^{pb}
\ee
Recall the identity (\ref{Nphi-switch}) for switching indices on
contractions of $N$ and $\varphi$. 
Therefore
\bea
N^\alpha{}_{kp} N_{\alpha}{}^{\lambda \ell}\varphi_{\lambda j \ell} \varphi_{iab}\o^{ka}\o^{pb} &=&  N^\alpha{}_{kp} [ N^\ell{}_{\lambda j} \varphi_{\ell \alpha}{}^\lambda + N^{\lambda}{}_{j \ell} \varphi_{\lambda}{}^\ell{}_\alpha ] \varphi_{iab}\o^{ka}\o^{pb} \nonumber\\
&=& -2  N^\alpha{}_{kp}  N^{\ell \lambda}{}_j \varphi_{\alpha \ell  \lambda} \varphi_{iab}\o^{ka}\o^{pb}.
\eea
Applying the identity (\ref{Nphi-switch}) again,
\be
N^\alpha{}_{kp} N_{\alpha}{}^{\lambda \ell}\varphi_{\lambda j \ell} \varphi_{iab}\o^{ka}\o^{pb} = 2  N^\alpha{}_{\ell \lambda}  N^{\ell \lambda}{}_j \varphi_{\alpha kp} \varphi_{iab}\o^{ka}\o^{pb}.
\ee
We can now apply the bilinear identities (\ref{contract}), so that
\be \label{F-terms3}
N^\alpha{}_{kp} N_{\alpha}{}^{\lambda \ell}\varphi_{\lambda j \ell} \varphi_{iab}\o^{ka}\o^{pb} =- 2 | \varphi |^2 N^\alpha{}_{\ell \lambda}  N^{\ell \lambda}{}_j g_{\alpha i} = - 2 | \varphi |^2 N_{i \ell \lambda}  N^{\ell \lambda}{}_j.
\ee
Next, we need to handle the ${\frak D} N$ terms in (\ref{F-terms}). By the Bianchi identity  $N_{ijk} + N_{kij} + N_{jki} = 0$, we have
\be
-2 {\frak D}_k N_{p}{}^{\lambda \ell} \varphi_{\lambda j \ell} \varphi_{iab}\o^{ka}\o^{pb} = 2 {\frak D}_k N^\ell{}_p{}^\lambda \varphi_{\lambda j \ell} \varphi_{iab}\o^{ka}\o^{pb} +2 {\frak D}_k N^{\lambda \ell}{}_p \varphi_{\lambda j \ell} \varphi_{iab}\o^{ka}\o^{pb}
\ee
This is
\be \label{F-terms2}
-2 {\frak D}_k N_{p}{}^{\lambda \ell} \varphi_{\lambda j \ell} \varphi_{iab}\o^{ka}\o^{pb} = -4 {\frak D}_k N^{\ell \lambda}{}_p \varphi_{\ell \lambda j} \varphi_{iab}\o^{ka}\o^{pb}.
\ee
To apply the bilinear identities (\ref{contract}), we will need to switch some indices.
\begin{lemma}
\be
{\frak D}_k N^p{}_{ij} \varphi_{p \lambda l}  = -{\frak D}_k N^p{}_{\lambda l} \varphi_{p ij} +  N^p{}_{\lambda, J l} \alpha_{Jk} \varphi_{p ij}  .
\ee
\end{lemma}
{\it Proof:} Differentiating identity (\ref{Nphi-switch}) gives
\be
{\frak D}_k N^p{}_{ij} \varphi_{p \lambda l} +  N^p{}_{ij} {\frak D}_k \varphi_{p \lambda l} = -{\frak D}_k N^p{}_{\lambda l} \varphi_{p ij} - N^p{}_{\lambda l}  {\frak D}_k \varphi_{p ij} .
\ee
Using the formula (\ref{frakD-varphi}), we obtain
\be
{\frak D}_k N^p{}_{ij} \varphi_{p \lambda l}  = -{\frak D}_k N^p{}_{\lambda l} \varphi_{p ij} + {1 \over 2}  N^p{}_{ij} \alpha_k \varphi_{p \lambda l} + {1 \over 2}  N^p{}_{ij} \alpha_{Jk} \hat{\varphi}_{p \lambda l} + {1 \over 2} N^p{}_{\lambda l}  \alpha_k \varphi_{p ij} + {1 \over 2} N^p{}_{\lambda l} \alpha_{Jk} \hat{\varphi}_{p ij} .
\ee
Using (\ref{Nphi-switch}) and $\hat{\varphi}_{p \lambda l} = -\varphi_{p J \lambda, l} = - \varphi_{Jp, \lambda l}$, we simplify this to
\be
{\frak D}_k N^p{}_{ij} \varphi_{p \lambda l}  = -{\frak D}_k N^p{}_{\lambda l} \varphi_{p ij} - {1 \over 2}  N^p{}_{ij} \alpha_{Jk} \varphi_{p, J \lambda, l}  - {1 \over 2} N^p{}_{\lambda l} \alpha_{Jk} \varphi_{J p, ij} .
\ee
Using (\ref{Nphi-switch}) again and $N^{Jp}{}_{\lambda l} = - N^p{}_{\lambda , Jl}$, we obtain the desired identity. Q.E.D.
\bigskip

Applying now this lemma to (\ref{F-terms2}), we find
\be
-2 {\frak D}_k N_{p}{}^{\lambda \ell} \varphi_{\lambda j \ell} \varphi_{iab}\o^{ka}\o^{pb} = (4 {\frak D}_k N^{\ell \lambda}{}_j  -  4 N^{\ell \lambda}{}_{Jj} \alpha_{Jk} ) \varphi_{\ell \lambda p} \varphi_{iab}\o^{ka}\o^{pb}
\ee
We can now use the bilinear identities
\bea
&&-2 {\frak D}_k N_{p}{}^{\lambda \ell} \varphi_{\lambda j \ell} \varphi_{iab}\o^{ka}\o^{pb} \nonumber\\
&=& | \varphi |^2 ( {\frak D}_k N^{\ell \lambda}{}_j - N^{\ell \lambda}{}_{Jj} \alpha_{Jk} ) (\omega_{\ell i} g_{\lambda a} -\omega_{\lambda i} g_{\ell a} -\omega_{\ell a} g_{\lambda i} + \omega_{\lambda a} g_{\ell i}  )\o^{ka} \nonumber\\
&=& | \varphi |^2 ( {\frak D}_k N^{\ell \lambda}{}_j - N^{\ell \lambda}{}_{Jj} \alpha_{Jk} ) (\omega_{\ell i} J^k{}_\lambda -\omega_{\lambda i} J^k{}_\ell +\delta^k{}_\ell g_{\lambda i} - \delta^k{}_\lambda g_{\ell i}  ) \nonumber\\
&=&  | \varphi |^2 (- {\frak D}_k N_{Ji}{}^{Jk}{}_j + N_{Ji}{}^{Jk}{}_{Jj} \alpha_{Jk} ) + | \varphi |^2 ( {\frak D}_k N^{Jk}{}_{Ji,j} - N^{Jk}{}_{Ji,Jj} \alpha_{Jk} ) \nonumber\\
&&+ | \varphi |^2 ( {\frak D}_k N^k{}_{ij} - N^{k}{}_{i,Jj} \alpha_{Jk} ) + | \varphi |^2 ( -{\frak D}_k N_i{}^{k}{}_j + N_i{}^{k}{}_{Jj} \alpha_{Jk} ).
\nonumber
\eea
This simplifies to
\be \label{F-terms4}
 -2 {\frak D}_k N_{p}{}^{\lambda \ell} \varphi_{\lambda j \ell} \varphi_{iab}\o^{ka}\o^{pb} = | \varphi |^2 (-2 {\frak D}_k N_i{}^k{}_j + 2 {\frak D}_k N^k{}_{ij} + 2 N_i{}^k{}_j \alpha_k - 2 N^k{}_{ij} \alpha_k).
\ee
Substituting (\ref{F-terms3}) and (\ref{F-terms4}) into (\ref{F-terms}),
\bea
&&(R^\ell{}_p{}^{\lambda}{}_k - R^\ell{}_k{}^{\lambda}{}_p)\varphi_{\lambda j \ell} \varphi_{iab}\o^{ka}\o^{pb} \nonumber\\
&=& | \varphi |^2 (-2 {\frak D}_k N_i{}^k{}_j + 2 {\frak D}_k N^k{}_{ij} + 2 N_i{}^k{}_j \alpha_k - 2 N^k{}_{ij} \alpha_k)
 + 2 | \varphi |^2 N_{i \ell \lambda}  N^{\ell \lambda}{}_j
\eea
By the Bianchi identity,
\be
2 | \varphi |^2 N_{i \ell \lambda}  N^{\ell \lambda}{}_j = 2 | \varphi |^2 ( - N_{\lambda i \ell} - N_{\ell \lambda i} )  N^{\ell \lambda}{}_j = 2 | \varphi |^2 N_{\lambda \ell i} N^{\ell \lambda}{}_j - 2 | \varphi |^2 N_{\ell \lambda i} N^{\ell \lambda}{}_j
\ee
and hence
\bea
(R^\ell{}_p{}^{\lambda}{}_k - R^\ell{}_k{}^{\lambda}{}_p)\varphi_{\lambda j \ell} \varphi_{iab}\o^{ka}\o^{pb} 
&=& | \varphi |^2 (-2 {\frak D}_k N_i{}^k{}_j + 2 {\frak D}_k N^k{}_{ij} + 2 N_i{}^k{}_j \alpha_k - 2 N^k{}_{ij} \alpha_k)
\nonumber\\
&& +  2 | \varphi |^2 (N^2_-)_{ij} - 2 | \varphi |^2 (N^2_+)_{ij} \nonumber
\eea
The result is
\bea
F &=&  | \varphi |^2 \bigg\{ (- R_{J i, J \lambda}{}^\lambda{}_j -
R_{J j, J \lambda}{}^\lambda{}_i) + (R_{Ji,j}{}^\lambda{}_{J\lambda} + R_{Jj,i}{}^\lambda{}_{J\lambda}) + 2 R_{ij} \nonumber\\
&& -2 ({\frak D}_k N_i{}^k{}_j + {\frak D}_k N_j{}^k{}_i) + 2 (N_i{}^k{}_j + N_j{}^k{}_i) \alpha_k + 4 (N^2_-)_{ij} - 4 (N^2_+)_{ij}   \bigg\}
\eea



\begin{lemma} \label{curv-contributions1}
We have the following formula
\bea
& \ & g^{\ell m}([\na_m,\na_j]\varphi_{kp\ell}
+
[\na_m,\na_k]\varphi_{pj\ell}+[\na_m,\na_p]\varphi_{jk\ell})
\varphi_{iab}\o^{ka}\o^{pb}+(i\leftrightarrow j)
\nonumber\\
&=& | \varphi |^2 \bigg\{ -2 R_{ji} -  R g_{ij} + R_{ij} + R_{Jj,Ji} \nonumber\\
&&-( R_{Ji, J \lambda}{}^\lambda{}_j + R_{Jj, J \lambda}{}^\lambda{}_i) + (R_{i,Jj}{}^\lambda{}_{J\lambda} + R_{j,Ji}{}^\lambda{}_{J\lambda}) + 2 R_{ij} \nonumber\\
&&  -2 ({\frak D}_k N_i{}^k{}_j + {\frak D}_k N_j{}^k{}_i) + 2 (N_i{}^k{}_j + N_j{}^k{}_i) \alpha_k   +4 (N^2_-)_{ij} - 4 (N^2_+)_{ij} \bigg\} 
\eea

\end{lemma}

\subsubsection{Contributions of the curvature terms, continued}
We now simplify Lemma \ref{curv-contributions1} by applying identities for the action of $J$ on the Riemann curvature tensor. We start with the terms
\be
- R_{Ji,J \lambda}{}^\lambda{}_j - R_{Jj,J \lambda}{}^\lambda{}_i
\ee
which can be manipulated using the relation (\ref{riemann-J-identity}) into
\bea \label{R-J1}
- R_{Ji,J \lambda}{}^\lambda{}_j - R_{Jj,J \lambda}{}^\lambda{}_i &=& - R_j{}^\lambda{}_{J \lambda, Ji} - R_i{}^\lambda{}_{J \lambda,Jj} \nonumber\\
&=& - R_j{}^\lambda{}_{\lambda i} - R_i{}^\lambda{}_{\lambda j} - B^\lambda{}_{j \lambda i} - B^\lambda{}_{i \lambda j} \nonumber\\
&=&  2 R_{ij}- B^\lambda{}_{j \lambda i} - B^\lambda{}_{i \lambda j}.
\eea
Next, we have the terms
\be
R_{i,Jj}{}^\lambda{}_{J \lambda} + R_{j,Ji}{}^\lambda{}_{J \lambda}.
\ee
By the Bianchi identity,
\bea
R_{i,Jj}{}^\lambda{}_{J \lambda} + (i \leftrightarrow j) &=& - R_{j,J
  \lambda}{}^\lambda{}_{J i} - R_{J\lambda, Ji}{}^\lambda{}_j + (i \leftrightarrow j) \nonumber\\
&=&  -  R_{j \lambda}{}^{J \lambda}{}_{J i} - R_j{}^\lambda{}_{Ji, J\lambda} + (i \leftrightarrow j) \nonumber\\
&=&  g^{\lambda \mu} R_{j, \lambda, J \mu, J i} - R_j{}^\lambda{}_{Ji, J\lambda} + (i \leftrightarrow j)
\eea
Using the relation (\ref{riemann-J-identity}),
\bea
R_{i,Jj}{}^\lambda{}_{J \lambda} + (i \leftrightarrow j) &=& g^{\lambda \mu} R_{j, \lambda,  \mu,  i} - R_j{}^\lambda{}_{i, \lambda} + g^{\lambda \mu} B_{ \lambda, j,  \mu,  i}  - B^\lambda{}_{j i \lambda} + (i \leftrightarrow j) \nonumber\\
&=& - 2 R_{ij} - 2 R_{ij} + \{  B^\lambda{}_{j \lambda i} - B^\lambda{}_{j i \lambda} +  B^\lambda{}_{i \lambda j} - B^\lambda{}_{i j \lambda}  \}
\eea
Therefore
\be \label{R-J2}
R_{Ji,j}{}^\lambda{}_{J \lambda} + R_{Jj,i}{}^\lambda{}_{J \lambda} = -4 R_{ij} + \{  B^\lambda{}_{j \lambda i} - B^\lambda{}_{j i \lambda}+  B^\lambda{}_{i \lambda j} - B^\lambda{}_{i j \lambda}  \}.
\ee
The next term in Lemma \ref{curv-contributions1} that we consider is $R_{Jj,Ji}$. This term becomes
\bea
R_{Jj,Ji} = g^{\lambda \mu} R_{\lambda, Jj, \mu, Ji} &=& - g^{\lambda \mu} R_{\lambda,Jj, J \mu, i} - g^{\lambda \mu} B_{Jj, \lambda, J \mu, i} \nonumber\\
&=& - g^{\lambda \mu} R_{i, J \mu, Jj, \lambda} - g^{\lambda \mu} B_{Jj, \lambda, J \mu, i} \nonumber\\
&=&   g^{\lambda \mu} R_{i, J \mu, j, J \lambda} + g^{\lambda \mu} B_{ J \mu, i, j, J\lambda}  - g^{\lambda \mu} B_{Jj, \lambda, J \mu, i}
\eea
and thus
\be \label{R-J3}
R_{Jj,Ji} = R_{ij} + B^\lambda{}_{ij \lambda} - B_{Jj}{}^\lambda{}_{J \lambda,i} .
\ee
Substituting (\ref{R-J1}), (\ref{R-J2}) and (\ref{R-J3}) into Lemma \ref{curv-contributions1}, we obtain
\bea
& \ & g^{\ell m}([\na_m,\na_j]\varphi_{kp\ell}
+
[\na_m,\na_k]\varphi_{pj\ell}+[\na_m,\na_p]\varphi_{jk\ell})
\varphi_{iab}\o^{ka}\o^{pb}+(i\leftrightarrow j)
\nonumber\\
&=& - | \varphi |^2 R \, g_{ij} -2 ({\frak D}_k N_i{}^k{}_j + {\frak D}_k N_j{}^k{}_i) + 2 (N_i{}^k{}_j + N_j{}^k{}_i) \alpha_k + 4  (N^2_-)_{ij} - 4 (N^2_+)_{ij}  \nonumber\\
&&  - B^\lambda{}_{j i \lambda}  - B_{Jj}{}^\lambda{}_{J \lambda,i}
\eea
Using the definition of $B$,
\bea
- B^\lambda{}_{ji \lambda} - B_{Jj}{}^\lambda{}_{J \lambda,i} &=& - [ - 2 D^\lambda N_{ji \lambda} + 2 D_j N^\lambda{}_{i \lambda} - 2N^{\alpha \lambda}{}_{j} N_{\alpha i \lambda}] \nonumber\\
&&- [ - 2 D_{Jj} N^\lambda{}_{J \lambda, i} + 2 D^\lambda N_{Jj,J\lambda,i} - 2 N^\alpha{}_{Jj}{}^\lambda N_{\alpha, J \lambda,i}] \nonumber\\
&=& - 4 (N^2_+)_{ij}
\eea
where we use the symmetries of $N$ to get the last equality. Therefore

\begin{lemma} \label{curv-contributions2}
We have the following formula
\bea
& \ & g^{\ell m}([\na_m,\na_j]\varphi_{kp\ell}
+
[\na_m,\na_k]\varphi_{pj\ell}+[\na_m,\na_p]\varphi_{jk\ell})
\varphi_{iab}\o^{ka}\o^{pb}+(i\leftrightarrow j)
\\
&=&  | \varphi |^2 \left\{ - R g_{ij} + 2 {\frak D}_k N_{ij}{}^k  + 2 {\frak D}_k N_{ji}{}^k + 2 (N_i{}^k{}_j + N_j{}^k{}_i) \alpha_k   
 + 4 (N^2_-)_{ij} - 8 (N^2_+)_{ij}\right\}
\nonumber
\eea

\end{lemma}

\subsubsection{Bochner-Kodaira contributions}
By (\ref{bk-identity1}), we have
\bea
 (- | \varphi |^2 d d^\dagger \varphi)_{jkp} \varphi_{iab}\o^{ka}\o^{pb}
&=& ( | \varphi |^2 g^{\ell m} \nabla_m \nabla_\ell \varphi_{jkp} ) \varphi_{iab}\o^{ka}\o^{pb} - | \varphi |^2 ( g^{\ell m} \big\{[\na_m,\na_j]\varphi_{kp\ell} \nonumber\\
&&\qquad+ [\na_m,\na_k]\varphi_{pj \ell} + [\na_m,\na_p]\varphi_{jk\ell}\big\}) \varphi_{iab}\o^{ka}\o^{pb}
\nonumber
\eea
By Lemma \ref{laplacian-contribution} and Lemma \ref{curv-contributions2}, we obtain
\bea 
&&
(- | \varphi |^2 d d^\dagger \varphi)_{jkp} \varphi_{iab}\o^{ka}\o^{pb} + (i \leftrightarrow j) 
= | \varphi|^4 \big\{ \nabla_\mu \alpha^\mu+ |N|^2\big\}g_{ij} \nonumber\\
&& +| \varphi |^4 \left\{ R g_{ij} + 2 (-{\frak D}_k N_{i j}{}^k  -  {\frak D}_k N_{ji}{}^k) - 2 (N_i{}^k{}_j + N_j{}^k{}_i) \alpha_k 
- 4 (N^2_-)_{ij} + 8 (N^2_+)_{ij} \right\}\nonumber
\eea
Altogether,
\begin{lemma}
We have the following formula
\bea \label{bochner-contribution}
&&(- | \varphi |^2 d d^\dagger \varphi)_{jkp} \varphi_{iab}\o^{ka}\o^{pb} + (i \leftrightarrow j) \\
&=& | \varphi |^4 \bigg\{ R g_{ij} - 2 ({\frak D}_k N_{i j}{}^k  +  {\frak D}_k N_{ji}{}^k) + ( \nabla_\mu \alpha^\mu+ |N|^2 ) g_{ij}
- 2 (N_i{}^k{}_j + N_j{}^k{}_i) \alpha_k \nonumber\\
&& - 4 (N^2_-)_{ij} + 8 (N^2_+)_{ij} \bigg\}.\nonumber
\eea
\end{lemma}

\subsection{Other contributions}
\subsubsection{Gradient dagger}

Returning to (\ref{phi-flow-rhs}), we study the contributions of the second term $-d|\varphi|^2\wedge d^\dagger\varphi$.
We let $\alpha = - d \log | \varphi |^2$ as before, and write
\be
- d | \varphi |^2 = | \varphi |^2 \alpha, \quad (d^\dagger \varphi)_{kp} = - g^{\mu \beta} \nabla_\beta \varphi_{\mu kp}.
\ee
Since
\be
(-d|\varphi|^2\wedge d^\dagger\varphi)_{jkp} = (- d | \varphi |^2)_j (d^\dagger \varphi)_{kp} + (- d | \varphi |^2)_p (d^\dagger \varphi)_{jk} + (- d | \varphi |^2)_k (d^\dagger \varphi)_{pj}
\ee
we have
\be
(-d|\varphi|^2\wedge d^\dagger\varphi)_{jkp} = | \varphi |^2 \left( - \alpha_j g^{\mu \beta} \nabla_\beta \varphi_{\mu kp} - \alpha_p  g^{\mu \beta} \nabla_\beta \varphi_{\mu jk} - \alpha_k g^{\mu \beta} \nabla_\beta \varphi_{\mu pj} \right)
\ee
Using previous notation,
\be
\nabla_\beta \varphi_{\mu kp} = {\frak D}_\beta \varphi_{\mu kp} - E_{\beta;\mu kp}.
\ee
By the formula (\ref{frakD-varphi}), we conclude
\be
\nabla_\beta \varphi_{\mu kp} = - {1 \over 2} \alpha_\beta \varphi_{\mu kp} + {1 \over 2} \alpha_{J \beta} \varphi_{J \mu, kp} - E_{\beta;\mu kp}.
\ee
Therefore
\bea
(-d|\varphi|^2\wedge d^\dagger\varphi)_{jkp} &=& | \varphi |^2 \bigg( {1 \over 2} \alpha_j g^{\mu \beta} \alpha_\beta \varphi_{\mu kp} - {1 \over 2}\alpha_j g^{\mu \beta} \alpha_{J \beta} \varphi_{J \mu, kp} + \alpha_j g^{\mu \beta} E_{\beta;\mu kp} \nonumber\\
&&+  {1 \over 2} \alpha_p g^{\mu \beta} \alpha_\beta \varphi_{\mu jk} - {1 \over 2}\alpha_p g^{\mu \beta} \alpha_{J \beta} \varphi_{J \mu, jk} + \alpha_p g^{\mu \beta} E_{\beta;\mu jk} \nonumber\\
&&+  {1 \over 2} \alpha_k g^{\mu \beta} \alpha_\beta \varphi_{\mu pj} - {1 \over 2}\alpha_k g^{\mu \beta} \alpha_{J \beta} \varphi_{J \mu, pj} + \alpha_k g^{\mu \beta} E_{\beta;\mu pj} \bigg)
\eea
which simplifies to
\bea
(-d|\varphi|^2\wedge d^\dagger\varphi)_{jkp} &=& | \varphi |^2 \left( \alpha_j g^{\mu \beta} E_{\beta;\mu kp} + \alpha_p g^{\mu \beta} E_{\beta;\mu jk} + \alpha_k g^{\mu \beta} E_{\beta;\mu pj} \right) \nonumber\\
&:=& {\rm (I)} + {\rm (II)} + {\rm (III)}.
\eea
We now work out the bilinears.
\be
( {\rm I} ) \cdot \varphi_{iab} \omega^{ka} \omega^{pb} =   | \varphi |^2 \alpha_j g^{\mu \beta} ( \varphi_{\lambda kp} N_{\beta \mu}{}^\lambda + \varphi_{\mu \lambda p} N_{\beta k}{}^\lambda + \varphi_{\mu k \lambda} N_{\beta p}{}^\lambda) \varphi_{iab} \omega^{ka} \omega^{pb} 
\ee
Since $N^\mu{}_{\mu}{}^\lambda = 0$ and we can relabel $p \leftrightarrow k$ and $a \leftrightarrow b$,
\be
( {\rm I}) \cdot \varphi_{iab} \omega^{ka} \omega^{pb} =  2  | \varphi |^2 \alpha_j g^{\mu \beta} ( \varphi_{\mu \lambda p} N_{\beta k}{}^\lambda) \varphi_{iab} \omega^{ka} \omega^{pb}.
\ee
By the bilinear identities
\bea
( {\rm I}) \cdot \varphi_{iab} \omega^{ka} \omega^{pb} &=& {| \varphi |^4 \over 2} \alpha_j g^{\mu \beta} N_{\beta k}{}^\lambda ( g_{\mu i} \omega_{\lambda a} - g_{\lambda i} \omega_{\mu a} - g_{\mu a} \omega_{\lambda i} + g_{\lambda a} \omega_{\mu i}) \omega^{ka}  \nonumber\\
&=& {| \varphi |^4 \over 2} \alpha_j g^{\mu \beta} N_{\beta k}{}^\lambda ( -g_{\mu i} \delta^k{}_\lambda + g_{\lambda i} \delta^k{}_\mu - J^k{}_\mu \omega_{\lambda i} + J^k{}_\lambda \omega_{\mu i}) \nonumber\\
&=& {| \varphi |^4 \over 2} ( 0+0 + \alpha_j N^{Jk}{}_k{}^{Ji} - \alpha_j N_{Ji,J\lambda}{}^\lambda) = 0
\eea
using the type $(0,2)$ and trace-free property of $N$. Next,
\bea
( {\rm II} + {\rm III} ) \cdot \varphi_{iab} \omega^{ka} \omega^{pb} &=&  2 | \varphi |^2 \alpha_p g^{\mu \beta} ( \varphi_{\lambda jk} N_{\beta \mu}{}^\lambda + \varphi_{\mu \lambda k} N_{\beta j}{}^\lambda + \varphi_{\mu j \lambda} N_{\beta k}{}^\lambda) \varphi_{iab} \omega^{ka} \omega^{pb} \nonumber\\
&=&  2 | \varphi |^2 \alpha_p  ( 0 + \varphi_{\mu \lambda k} N^\mu{}_{j}{}^\lambda + \varphi_{\mu j \lambda} N^\mu{}_{k}{}^\lambda) \varphi_{iab} \omega^{ka} \omega^{pb}
\eea
The first term is
\bea
2 | \varphi |^2 \alpha_p  (\varphi_{\mu \lambda k} N^\mu{}_{j}{}^\lambda) \varphi_{iab} \omega^{ka} \omega^{pb} &=& - 2 | \varphi |^2 \alpha_p N^\mu{}_{j}{}^\lambda (\varphi_{\mu \lambda k}  \varphi_{iba} \omega^{ka}) \omega^{pb}\nonumber\\
&=& -{| \varphi |^4 \over 2} \alpha_p N^\mu{}_{j}{}^\lambda ( g_{\mu i} \omega_{\lambda b} - g_{\lambda i} \omega_{\mu b} - g_{\mu b} \omega_{\lambda i} + g_{\lambda b} \omega_{\mu i}) \omega^{pb} \nonumber\\
&=& -{| \varphi |^4 \over 2} \alpha_p N^\mu{}_{j}{}^\lambda ( - g_{\mu i} \delta^p{}_\lambda + g_{\lambda i} \delta^p{}_\mu - J^p{}_\mu \omega_{\lambda i} + J^p{}_\lambda \omega_{\mu i}) \nonumber\\
&=& -{| \varphi |^4 \over 2} ( - \alpha_p N_{ij}{}^p + \alpha_p N^p{}_{j i} +  \alpha_p N^{Jp}{}_{j, Ji} - \alpha_p N_{Ji, j}{}^{Jp} ) \nonumber\\
&=& | \varphi |^4 ( \alpha_p N_{ij}{}^p - \alpha_p N^p{}_{j i})
\eea
For the second term, we use the identity (\ref{Nphi-switch}) to obtain
\bea
2 | \varphi |^2 \alpha_p  (\varphi_{\mu j \lambda} N^\mu{}_{k}{}^\lambda) \varphi_{iab} \omega^{ka} \omega^{pb}&=& 2 | \varphi |^2 \alpha_p  ( -\varphi_{\mu k \lambda} N^\mu{}_{j}{}^\lambda ) \varphi_{iab} \omega^{ka} \omega^{pb} \nonumber\\
&=& - 2 | \varphi |^2 \alpha_p  N^\mu{}_{j}{}^\lambda  \varphi_{\mu \lambda k} \varphi_{iba} \omega^{ka}\omega^{pb} .
\eea
This term is identical to the one above. Therefore
\be
( {\rm II} + {\rm III} ) \cdot \varphi_{iab} \omega^{ka} \omega^{pb} = 2 | \varphi |^4 ( \alpha_p N_{ij}{}^p - \alpha_p N^p{}_{j i})
\ee
Altogether,
\be
(-d|\varphi|^2\wedge d^\dagger\varphi)_{jkp} \varphi_{iab} \omega^{ka} \omega^{pb} =  2 | \varphi |^4 (  \alpha_p N_{ij}{}^p - \alpha_p N^p{}_{j i}).
\ee
Therefore
  \be
(-d|\varphi|^2\wedge d^\dagger\varphi)_{jkp} \varphi_{iab} \omega^{ka} \omega^{pb} + (i \leftrightarrow j) =  2 | \varphi |^4 (  \alpha_p N_{ij}{}^p - \alpha_p N^p{}_{j i} + \alpha_p N_{ji}{}^p - \alpha_p N^p{}_{ij}).
  \ee
  By the Bianchi identity $ N^p{}_{ij} + N_j{}^p{}_i + N_{ij}{}^p = 0$, and hence
 \bea
&& (-d|\varphi|^2\wedge d^\dagger\varphi)_{jkp} \varphi_{iab} \omega^{ka} \omega^{pb} + (i \leftrightarrow j)\nonumber\\
 &=&  2 | \varphi |^4 \left(  \alpha_p N_{ij}{}^p - \alpha_p (-N_i{}^p{}_j - N_{ji}{}^p) + \alpha_p N_{ji}{}^p - \alpha_p (- N_j{}^p{}_i - N_{ij}{}^p)\right) \nonumber\\
 &=& 2| \varphi |^4 \alpha_p ( N_{ij}{}^p  + N_i{}^p{}_j + N_{ji}{}^p + N_{ji}{}^p + N_j{}^p{}_i + N_{ij}{}^p)
\eea
Thus we have established the following lemma:

\begin{lemma}
  \be\label{grad-dagger}
(-d|\varphi|^2\wedge d^\dagger\varphi)_{jkp} \varphi_{iab} \omega^{ka} \omega^{pb} + (i \leftrightarrow j) = - 2 | \varphi |^4 \alpha_p(N_{j}{}^p{}_i + N_{i}{}^p{}_j ).
  \ee
\end{lemma}

\subsubsection{Interior product}

Returning to (\ref{phi-flow-rhs}), we study the contributions of the third term
$d ( \iota_{\nabla | \varphi |^2} \varphi)$. We can write
\be
( \iota_{\nabla | \varphi |^2} \varphi)_{kp} = g^{\mu \nu} (\partial_\nu | \varphi |^2) \varphi_{\mu kp} = - | \varphi |^2 g^{\mu \nu} \alpha_\nu \varphi_{\mu kp}
\ee
since $\alpha_i = - \partial_i \log | \varphi |^2$, or $\partial_j | \varphi |^2 = - | \varphi |^2 \alpha_j$. Next,
\be \label{d-iota1}
d ( \iota_{\nabla | \varphi |^2} \varphi)_{jkp} = \nabla_j ( \iota_{\nabla | \varphi |^2} \varphi)_{kp} + \nabla_p ( \iota_{\nabla | \varphi |^2} \varphi)_{jk} + \nabla_k ( \iota_{\nabla | \varphi |^2} \varphi)_{pj}.
\ee
We start with
\be
\nabla_j ( \iota_{\nabla | \varphi |^2} \varphi)_{kp} = | \varphi |^2 \alpha_j \alpha^\mu \varphi_{\mu kp}- | \varphi |^2  \nabla_j \alpha^\mu \varphi_{\mu kp}- | \varphi |^2 \alpha^\mu \nabla_j \varphi_{\mu kp}
\ee
Since
\be
\nabla_j \varphi_{\mu kp} = - {1 \over 2} \alpha_j \varphi_{\mu kp} + {1 \over 2} \alpha_{J j} \varphi_{J \mu, k, p} - E_{j;\mu k p},
\ee
we have
\be
\nabla_j ( \iota_{\nabla | \varphi |^2} \varphi)_{kp} = {3 \over 2} | \varphi |^2 \alpha_j \alpha^\mu \varphi_{\mu kp}- | \varphi |^2  \nabla_j \alpha^\mu \varphi_{\mu kp} -  {1 \over 2} | \varphi |^2 \alpha^\mu  \alpha_{J j} \varphi_{J \mu, k, p}  + | \varphi |^2 \alpha^\mu   E_{j;\mu k p}.
\ee
We now work out the bilinears.
\be
({3 \over 2} | \varphi |^2 \alpha_j \alpha^\mu \varphi_{\mu kp} ) \varphi_{iab} \omega^{ka} \omega^{pb} = - {3 \over 2} | \varphi |^4 \alpha_j \alpha^\mu g_{\mu i} = - {3 \over 2} | \varphi |^4 \alpha_i \alpha_j,
\ee
\be
(- | \varphi |^2  \nabla_j \alpha^\mu \varphi_{\mu kp} ) \varphi_{iab} \omega^{ka} \omega^{pb} = | \varphi |^4 \nabla_j \alpha_i,
\ee
\be
(-  {1 \over 2} | \varphi |^2 \alpha^\mu  \alpha_{J j} \varphi_{J \mu, k, p}) \varphi_{iab} \omega^{ka} \omega^{pb} = {1 \over 2} | \varphi |^4 \alpha^\mu  \alpha_{J j} g_{J \mu, i} = - {1 \over 2} | \varphi |^4 \alpha_{Ji} \alpha_{Jj}.
\ee
Therefore
\be \label{d-iota2}
(\nabla_j ( \iota_{\nabla | \varphi |^2} \varphi)_{kp})\varphi_{iab} \omega^{ka} \omega^{pb} = - {3 \over 2} | \varphi |^4 \alpha_i \alpha_j + | \varphi |^4 \nabla_j \alpha_i - {1 \over 2} | \varphi |^4 \alpha_{Ji} \alpha_{Jj} + | \varphi |^2 \alpha^\mu   E_{j;\mu k p}\varphi_{iab} \omega^{ka} \omega^{pb}.
\ee
Next, we work out the two next contributions of this term with the indices $(jkp)$ cyclically permuted. After forming bilinears, these two extra terms are identical.
\be
(\nabla_p ( \iota_{\nabla | \varphi |^2} \varphi)_{jk})\varphi_{iab} \omega^{ka} \omega^{pb} + (\nabla_k ( \iota_{\nabla | \varphi |^2} \varphi)_{pj})\varphi_{iab} \omega^{ka} \omega^{pb} = 2 (\nabla_p ( \iota_{\nabla | \varphi |^2} \varphi)_{jk})\varphi_{iab} \omega^{ka} \omega^{pb}
\ee
As before, we have
\be
\nabla_p ( \iota_{\nabla | \varphi |^2} \varphi)_{jk} =  {3 \over 2} | \varphi |^2 \alpha_p \alpha^\mu \varphi_{\mu jk}- | \varphi |^2  \nabla_p \alpha^\mu \varphi_{\mu jk} -  {1 \over 2} | \varphi |^2 \alpha^\mu  \alpha_{J p} \varphi_{J \mu, j, k}  + | \varphi |^2 \alpha^\mu   E_{p;\mu jk}
\ee
Forming bilinears,
\bea
({3 \over 2} | \varphi |^2 \alpha_p \alpha^\mu \varphi_{\mu jk}) \varphi_{iab} \omega^{ka} \omega^{pb} &=& -{3 \over 8} | \varphi |^4 \alpha_p \alpha^\mu (\omega_{\mu i} g_{jb} - \omega_{j i} g_{\mu b} - \omega_{\mu b} g_{ji} + \omega_{jb} g_{\mu i}) \omega^{pb} \nonumber\\
&=& {3 \over 8} | \varphi |^4 \alpha_p \alpha^\mu ( -\omega_{\mu i} J^p{}_j + \omega_{j i} J^p{}_\mu - \delta^p{}_\mu g_{ji} + \delta^p{}_j g_{\mu i}) \nonumber\\
&=& {3 \over 8} | \varphi |^4 ( \alpha_{Jj} \alpha_{Ji} + \alpha_{J \mu} \alpha^\mu \omega_{ij} - \alpha_\mu \alpha^\mu g_{ij} + \alpha_j \alpha_i),
\eea
and
\bea
(- | \varphi |^2  \nabla_p \alpha^\mu \varphi_{\mu jk}) \varphi_{iab} \omega^{ka} \omega^{pb} &=& {1 \over 4} | \varphi |^4 \nabla_p \alpha^\mu (\omega_{\mu i} g_{jb} - \omega_{j i} g_{\mu b} - \omega_{\mu b} g_{ji} + \omega_{jb} g_{\mu i}) \omega^{pb} \nonumber\\
&=& {1 \over 4} | \varphi |^4 \nabla_p \alpha^\mu ( \omega_{\mu i} J^p{}_j - \omega_{j i} J^p{}_\mu + \delta^p{}_\mu g_{ji} - \delta^p{}_j g_{\mu i}) \nonumber\\
&=& { | \varphi |^4 \over 4} (- J^n{}_j \nabla_{n} \alpha_{q} J^q{}_i - J^p{}_\mu \nabla_{p} \alpha^\mu \omega_{ji} + \nabla_\mu \alpha^\mu g_{ij} - \nabla_j \alpha_i), \nonumber
\eea
and
\bea
(-{1 \over 2} | \varphi |^2  \alpha^\mu \alpha_{Jp} \varphi_{J \mu, j, k}) \varphi_{iab} \omega^{ka} \omega^{pb} &=& {1 \over 8} | \varphi |^4 \alpha^\mu \alpha_{Jp} (\omega_{J \mu, i} g_{jb} - \omega_{j i} g_{J \mu, b} - \omega_{J \mu, b} g_{ji} + \omega_{jb} g_{J \mu, i}) \omega^{pb} \nonumber\\
&=& {1 \over 8} | \varphi |^4 \alpha^\mu \alpha_{Jp} ( - g_{\mu i} J^p{}_j + \omega_{j i} \delta^p{}_\mu + J^p{}_\mu g_{ji} - \delta^p{}_j g_{J \mu, i}) \nonumber\\
&=& { | \varphi |^4 \over 8} (\alpha_i \alpha_j + \alpha^p \alpha_{Jp} \omega_{ji} - \alpha^\mu \alpha_\mu g_{ij} + \alpha_{Ji} \alpha_{Jj}).
\eea
Altogether,
\bea
(\nabla_p ( \iota_{\nabla | \varphi |^2} \varphi)_{jk})  \varphi_{iab} \omega^{ka} \omega^{pb} &=& {| \varphi |^4 \over 8} \bigg(  3 \alpha_{Jj} \alpha_{Ji} + 3 \alpha_{J \mu} \alpha^\mu \omega_{ij} - 3 \alpha_\mu \alpha^\mu g_{ij} + 3 \alpha_j \alpha_i \nonumber\\
&&- 2 \nabla_{Jj} \alpha_{Ji} - 2 \nabla_{J \mu} \alpha^\mu \omega_{ji} + 2 \nabla_\mu \alpha^\mu g_{ij} - 2 \nabla_j \alpha_i \nonumber\\
&& + \alpha_i \alpha_j + \alpha^p \alpha_{Jp} \omega_{ji} - \alpha^\mu \alpha_\mu g_{ij} + \alpha_{Ji} \alpha_{Jj} \bigg) \nonumber\\
&&+  | \varphi |^2 \alpha^\mu   E_{p;\mu jk} \varphi_{iab} \omega^{ka} \omega^{pb}
\eea
It follows that
\bea \label{d-iota3}
2 (\nabla_p ( \iota_{\nabla | \varphi |^2} \varphi)_{jk})  \varphi_{iab} \omega^{ka} \omega^{pb} &=& {| \varphi |^4 \over 4} \bigg(  4 \alpha_{Jj} \alpha_{Ji} + 2 \alpha_{J \mu} \alpha^\mu \omega_{ij} - 4 \alpha_\mu \alpha^\mu g_{ij} + 4 \alpha_j \alpha_i \nonumber\\
&&- 2 J^n{}_j \nabla_{n} \alpha_{q} J^q{}_i - 2 J^p{}_\mu \nabla_{p} \alpha^\mu \omega_{ji} + 2 \nabla_\mu \alpha^\mu g_{ij} - 2 \nabla_j \alpha_i \bigg) \nonumber\\
&&+  2 | \varphi |^2 \alpha^\mu   E_{p;\mu jk} \varphi_{iab} \omega^{ka} \omega^{pb}
\eea
We can now combine all of our calculations. By (\ref{d-iota1}), (\ref{d-iota2}), (\ref{d-iota3}),

\begin{lemma} 
\bea
& \ & (d \iota_{\nabla | \varphi |^2} \varphi)_{jkp}  \varphi_{iab} \omega^{ka} \omega^{pb} + (i \leftrightarrow j) \nonumber\\
&=& | \varphi |^4 \big\{ {1 \over 2}( \nabla_j \alpha_i + \nabla_i \alpha_j)- \alpha_i \alpha_j  + \alpha_{Ji} \alpha_{Jj}  -  2 \alpha_\mu \alpha^\mu g_{ij}  - {1 \over 2} (J^p{}_j J^q{}_i \nabla_{p} \alpha_{q} + J^p{}_i J^q{}_j \nabla_{p} \alpha_{q}) \nonumber\\
&&+  \nabla_\mu \alpha^\mu g_{ij} \big\} + | \varphi |^2 \alpha^\mu   E_{j;\mu k p}\varphi_{iab} \omega^{ka} \omega^{pb} +  | \varphi |^2 \alpha^\mu   E_{i;\mu k p}\varphi_{jab} \omega^{ka} \omega^{pb} \nonumber\\
&& +  2| \varphi |^2 \alpha^\mu   E_{p;\mu jk} \varphi_{iab} \omega^{ka} \omega^{pb}+  2| \varphi |^2 \alpha^\mu   E_{p;\mu ik} \varphi_{jab} \omega^{ka} \omega^{pb}.
\eea
\end{lemma}

It remains to evaluate the $E$ terms.
\be
 | \varphi |^2 \alpha^\mu   E_{j;\mu k p}\varphi_{iab} \omega^{ka} \omega^{pb} + 2| \varphi |^2 \alpha^\mu   E_{p;\mu jk} \varphi_{iab} \omega^{ka} \omega^{pb} + (i \leftrightarrow j)
\ee
We start with
\be
 | \varphi |^2 \alpha^\mu   E_{j;\mu k p}\varphi_{iab} \omega^{ka} \omega^{pb} =  | \varphi |^2 \alpha^\mu  ( \varphi_{\lambda kp} N_{j \mu }{}^\lambda + \varphi_{\mu \lambda p} N_{j k}{}^\lambda + \varphi_{\mu k \lambda} N_{j p}{}^\lambda )  \varphi_{iab} \omega^{ka} \omega^{pb}
\ee
which by symmetry is
\be
 | \varphi |^2 \alpha^\mu   E_{j;\mu k p}\varphi_{iab} \omega^{ka} \omega^{pb} =  | \varphi |^2 \alpha^\mu  ( \varphi_{\lambda kp} N_{j \mu }{}^\lambda) \varphi_{iab} \omega^{ka} \omega^{pb} +  2  | \varphi |^2 \alpha^\mu(\varphi_{\mu \lambda p} N_{j k}{}^\lambda)  \varphi_{iab} \omega^{ka} \omega^{pb}
\ee
The first term is
\bea
| \varphi |^2 \alpha^\mu  ( \varphi_{\lambda kp} N_{j \mu }{}^\lambda) \varphi_{iab} \omega^{ka} \omega^{pb} &=& -| \varphi |^4 \alpha^\mu  N_{j \mu }{}^\lambda g_{\lambda i} 
= - | \varphi |^4 \alpha^\mu N_{j \mu i}.
\eea
The second term is
\bea
2 | \varphi |^2 \alpha^\mu (\varphi_{\mu \lambda p} N_{j k}{}^\lambda)  \varphi_{iab} \omega^{ka} \omega^{pb} &=&  {| \varphi |^2 \over 4} \alpha^\mu N_{j k}{}^\lambda (g_{\mu i} \omega_{\lambda a} - g_{\lambda i} \omega_{\mu a} - g_{\mu a} \omega_{\lambda i} + g_{\lambda a} \omega_{\mu i}) \omega^{ka} \nonumber\\
&=&  {| \varphi |^2 \over 4} \alpha^\mu N_{j k}{}^\lambda (- g_{\mu i} \delta^k{}_\lambda + g_{\lambda i} \delta^k{}_\mu - J^k{}_\mu \omega_{\lambda i} + J^k{}_\mu \omega_{\mu i}) \nonumber\\
&=&  {| \varphi |^2 \over 4} \alpha^\mu N_{j k}{}^\lambda (- g_{\mu i} \delta^k{}_\lambda + g_{\lambda i} \delta^k{}_\mu - J^k{}_\mu \omega_{\lambda i} + J^k{}_\lambda \omega_{\mu i}) \nonumber\\
&=& {| \varphi |^2 \over 4}   (- \alpha_i N_{j \lambda}{}^\lambda  + \alpha^\mu N_{j \mu i}  + \alpha^\mu N_{j , J \mu, Ji} - \alpha_{Ji} N_{j , J\lambda }{}^\lambda ) \nonumber\\
&=& {| \varphi |^2 \over 4}   ( 0 + \alpha^\mu N_{j \mu i}  - \alpha^\mu N_{j  \mu i} + 0 ) = 0.
\eea
Therefore
\be \label{innerproductE1}
| \varphi |^2 \alpha^\mu   E_{j;\mu k p}\varphi_{iab} \omega^{ka} \omega^{pb} = - | \varphi |^4 \alpha^\mu N_{j \mu i} = | \varphi |^4 \alpha^\mu N_{ji \mu} .
\ee
Next, we consider
\bea \label{innerproductE2}
2| \varphi |^2 \alpha^\mu   E_{p;\mu jk} \varphi_{iab} \omega^{ka} \omega^{pb} &=& 2| \varphi |^2 \alpha^\mu (\varphi_{\lambda jk} N_{p \mu}{}^\lambda + \varphi_{\mu \lambda k} N_{p j}{}^\lambda + \varphi_{\mu j \lambda} N_{p k}{}^\lambda ) \varphi_{iab} \omega^{ka} \omega^{pb} \nonumber\\
&:=& (\tilde {\rm I} +\tilde {\rm II} + \tilde{\rm III} ) \cdot \varphi_{iab} \omega^{ka} \omega^{pb}.
\eea
We start with
\bea
(\tilde {\rm I} ) \cdot \varphi_{iab} \omega^{ka} \omega^{pb} &=& -2| \varphi |^2 \alpha^\mu N_{p \mu}{}^\lambda (\varphi_{\lambda kj}  \varphi_{iab} \omega^{ka}) \omega^{pb} \nonumber\\
&=& -{| \varphi |^4 \over 2} \alpha^\mu N_{p \mu}{}^\lambda ( \omega_{\lambda i} g_{jb} - \omega_{ji} g_{\lambda b} - \omega_{\lambda b} g_{ji} + \omega_{jb} g_{\lambda i} ) \omega^{pb} \nonumber\\
&=& -{| \varphi |^4 \over 2} \alpha^\mu N_{p \mu}{}^\lambda ( \omega_{\lambda i} J^p{}_j - \omega_{ji} J^p{}_\lambda + \delta^p{}_\lambda g_{ji} - \delta^p{}_j g_{\lambda i} ) \nonumber\\
&=&  {| \varphi |^4 \over 2}  ( \alpha^\mu N_{Jj, \mu, Ji} + \omega_{ji} \alpha^\mu N_{J \lambda, \mu}{}^\lambda -  g_{ji} \alpha^\mu N_{\lambda \mu}{}^\lambda + \alpha^\mu N_{j \mu i} ) \nonumber\\
&=&  {| \varphi |^4 \over 2}  ( - \alpha^\mu N_{j \mu i} + 0 - 0 + \alpha^\mu N_{j \mu i} ) = 0.
\eea
Similarly, we can also compute
\bea
(\tilde {\rm II} ) \cdot \varphi_{iab} \omega^{ka} \omega^{pb} =0
\eea
The third term is
\be
(\tilde {\rm III} ) \cdot \varphi_{iab} \omega^{ka} \omega^{pb} =  2| \varphi |^2 \alpha^\mu (\varphi_{\mu j \lambda} N_{p k}{}^\lambda ) \varphi_{iab} \omega^{ka} \omega^{pb}
\ee
It can be rearranged using the symmetry $p \leftrightarrow k$, $a \leftrightarrow b$
\bea
( \tilde {\rm III} ) \cdot \varphi_{iab} \omega^{ka} \omega^{pb} &=&  2| \varphi |^2 \alpha^\mu \varphi_{\mu j \lambda} {(N_{p k}{}^\lambda - N_{k p}{}^\lambda )\over 2}  \varphi_{iab} \omega^{ka} \omega^{pb} \nonumber\\
&=& - | \varphi |^2 \alpha^\mu \varphi_{\lambda \mu j} N^\lambda{}_{pk} \varphi_{iab} \omega^{ka} \omega^{pb}
\eea
and then using the Bianchi identity. By the identity $N^\lambda{}_{pk} \varphi_{\lambda \mu j} = -N^\lambda{}_{\mu j} \varphi_{\lambda p k} $, 
\be
( \tilde {\rm III} ) \cdot \varphi_{iab} \omega^{ka} \omega^{pb} =  - | \varphi |^2 \alpha^\mu N^\lambda{}_{\mu j} \varphi_{\lambda p k}  \varphi_{iba} \omega^{ka} \omega^{pb}.
\ee
We can now use the bilinear identity.
\be
( \tilde {\rm III} ) \cdot \varphi_{iab} \omega^{ka} \omega^{pb} =  | \varphi |^4  \alpha^\mu N^\lambda{}_{\mu j} g_{\lambda i} = | \varphi |^4 \alpha^\mu N_{i \mu j} = - | \varphi|^4 \alpha^\mu N_{i  j \mu} .
\ee
Substituting our results into (\ref{innerproductE2}), we obtain
\be
2| \varphi |^2 \alpha^\mu   E_{p;\mu jk} \varphi_{iab} \omega^{ka} \omega^{pb} = - | \varphi |^4 \alpha^\mu N_{i j \mu} .
\ee
Combining the above equation with (\ref{innerproductE1}),
\bea
& \ &  | \varphi |^2 \alpha^\mu   E_{j;\mu k p}\varphi_{iab} \omega^{ka} \omega^{pb} + 2| \varphi |^2 \alpha^\mu   E_{p;\mu jk} \varphi_{iab} \omega^{ka} \omega^{pb} + (i \leftrightarrow j) \nonumber\\
&=&  | \varphi |^4 \alpha^\mu N_{j i \mu} - | \varphi |^4 \alpha^\mu N_{i j \mu} + (i \leftrightarrow j)
= 0.
\eea
Therefore the $E$ terms do not contribute, and we are left with:

\begin{lemma} \label{d-iota-contribution} 
\bea
(d \iota_{\nabla | \varphi |^2} \varphi)_{jkp}  \varphi_{iab} \omega^{ka} \omega^{pb} + (i \leftrightarrow j) &=&
 | \varphi |^4 \big\{ {1 \over 2}( \nabla_j \alpha_i + \nabla_i \alpha_j)- \alpha_i \alpha_j 
 + \alpha_{Ji} \alpha_{Jj}  -  2 \alpha_\mu \alpha^\mu g_{ij} \nonumber\\
 &&
 - {1 \over 2} (J^p{}_j J^q{}_i \nabla_{p} \alpha_{q} + J^p{}_i J^q{}_j \nabla_{p} \alpha_{q})+  \nabla_\mu \alpha^\mu g_{ij} \big\}
 \nonumber
\eea
\end{lemma}

\medskip

\subsection{$N^\dagger$ term: $d(|\varphi|^2N^\dagger\cdot\varphi)$}

Recall from the definition of the operator $N^\dagger$ that $(N^\dagger\varphi)_{kj}=2N^{\mu}{}_j{}^\lambda \varphi_{\mu k\lambda} $ , and thus
\bea\label{Ndagger}
d(|\varphi|^2N^\dagger\cdot\varphi)_{jkp}
&=&\nabla_j (|\varphi|^2(N^\dagger\cdot\varphi)_{kp}) + \nabla_p (|\varphi|^2(N^\dagger\cdot\varphi)_{jk}) +\nabla_k (|\varphi|^2(N^\dagger\cdot\varphi)_{pj}) \nonumber\\
&:=& {\rm I}+ {\rm II}+{\rm III}.
\eea

\subsubsection{Computation for (I)}
We start with the first term
\bea\label{Ndagger1}
\nabla_j (|\varphi|^2(N^\dagger\cdot\varphi)_{kp}) 
&=&
 -2 |\varphi|^2 \, \alpha_j N^\mu{}_p{}^\lambda \varphi_{\mu k \lambda} +2 |\varphi|^2 \nabla_j (N^\mu{}_p{}^\lambda \varphi_{\mu k \lambda}) \\
 &=&
 -2 |\varphi|^2 \, \alpha_j N^\mu{}_p{}^\lambda \varphi_{\mu k \lambda} +2 |\varphi|^2 \nabla_j N^\mu{}_p{}^\lambda \varphi_{\mu k \lambda} 
 +
 2 |\varphi|^2 N^\mu{}_p{}^\lambda \nabla_j \varphi_{\mu k \lambda}\nonumber
\eea
We now work out the bilinears term by term
\bea
-2|\varphi|^2 \, \alpha_j N^\mu{}_k{}^\lambda \varphi_{\mu p \lambda} \varphi_{iab}\o^{ka} \o^{pb}
&=&
{1\over 2}|\varphi|^4 \, \alpha_j N^\mu{}_k{}^\lambda \, ( \o_{\mu i} g_{\lambda a} + \o_{\lambda a} g_{\mu i} - \o_{\mu a} g_{\lambda i} - \o_{\lambda i} g_{\mu a} ) \o^{ka}\nonumber\\
&=&
{1\over 2}|\varphi|^4 \, \alpha_j N^\mu{}_k{}^\lambda \, (\o_{\mu i} J^k{}_\lambda - \delta^k{}_\lambda g_{\mu i} + \delta^k{}_\mu g_{\lambda i} - \o_{\lambda i} J^k{}_\mu) 
\nonumber\\
&=& 
{1\over 2} |\varphi|^4 \, \alpha_j (- N_{Ji, k}{}^{Jk} - N_{i k}{}^k + N^k{}_{k i} + N^{Jk}{}_{k, Ji})=0.\nonumber
\eea
The first two terms are zero due to anti-symmetry of $N$ in the second and third indices. The third and fourth terms are also zero since $g^{ml}N_{mlj} =0$ and $N^{Jk}{}_{k, Ji} = - N^k{}_{Jk, Ji} = N^k{}_{ki}$.

Next, we work with the second group of terms in (\ref{Ndagger1}):
\bea
2 |\varphi|^2 \nabla_j N^\mu{}_p{}^\lambda \varphi_{\mu k \lambda} \varphi_{iab} \o^{ka} \o^{pb}
&=&
{1\over 2} |\varphi|^4\nabla_j N^\mu{}_p{}^\lambda  ( \o_{\mu i} g_{\lambda b} + \o_{\lambda b} g_{\mu i} - \o_{\mu b} g_{\lambda i} - \o_{\lambda i} g_{\mu b}) \o^{pb}\nonumber\\
&=&
{1\over 2} |\varphi|^4 \nabla_j N^\mu{}_p{}^\lambda ( \o_{\mu i} J^p{}_\lambda - \delta^p{}_\lambda g_{\mu i} + \delta^p{}_\mu g_{\lambda i} - \o_{\lambda i} J^p{}_\mu)\nonumber\\
&=& {1\over 2} |\varphi|^4 (\nabla_j N^\mu{}_p{}^\lambda  \o_{\mu i} J^p{}_\lambda  - \nabla_j N_{i p}{}^p + \nabla_j N^p{}_{p i} - \nabla_j N^\mu{}_p{}^\lambda \o_{\lambda i} J^p{}_\mu)\nonumber\\
&=&
{1\over 2} |\varphi|^4 (\nabla_j N^\mu{}_p{}^\lambda  \o_{\mu i} J^p{}_\lambda - \nabla_j N^\mu{}_p{}^\lambda \o_{\lambda i} J^p{}_\mu)
\eea
The last two terms require extra work since $J$ may not be covariantly constant under $\nabla$.
\bea
 \o_{\mu i} \nabla_j N^\mu{}_p{}^\lambda J^p{}_\lambda &=&
 \o_{\mu i}( \nabla_j (N^\mu{}_{p}{}^\lambda J^p{}_\lambda ) - N^{\mu}{}_p{}^\lambda \nabla_j J^p{}_\lambda)
 =
 \o_{\mu i} (\nabla_j N^\mu{}_p{}^{Jp} +2 N^\mu{}_p{}^\lambda \, N_{j\lambda}{}^{Jp})\nonumber\\
 &=& 
 2 \o_{\mu i} N^\mu{}_p{}^\lambda \, N_{j\lambda}{}^{Jp}= - 2 N_{Ji, p}{}^\lambda \, N_{j\lambda}{}^{Jp}
 =- 2 N_{i, Jp}{}^\lambda \, N_{j\lambda}{}^{Jp}
 = 2 N_{ip}{}^\lambda N_{j\lambda}{}^p.
 \nonumber
\eea
Similarly, we can compute
\bea
 - \nabla_j N^\mu{}_p{}^\lambda \o_{\lambda i} J^p{}_\mu=  2 N^\mu{}_{Jp, i} \, N_{j \mu}{}^{Jp}=-2N^\mu{}_{pi} \, N_{j \mu}{}^{p}.
\eea
Altogether, we have
\bea
2|\varphi|^2 \nabla_j N^\mu{}_p{}^\lambda \varphi_{\mu k \lambda}  \varphi_{iab} \o^{ka} \o^{pb}
&=& 
|\varphi|^4(N_{ip}{}^\lambda N_{j\lambda}{}^p-N^\lambda{}_{pi} \, N_{j \lambda}{}^{p}).
\eea
Next, we consider the last group of terms in (\ref{Ndagger1}). 
\bea\label{3rd-term}
 2|\varphi|^2  N^\mu{}_p{}^\lambda\nabla_j \varphi_{\mu k \lambda} \varphi_{iab} \o^{ka} \o^{pb}
\eea
Since
\be
\nabla_j \varphi_{\mu k\lambda} = - {1 \over 2} \alpha_j \varphi_{\mu k\lambda} + {1 \over 2} \alpha_{J j} \varphi_{J \mu, k\lambda} - E_{j;\mu k \lambda},
\ee
then
\bea
2|\varphi|^2  N^\mu{}_p{}^\lambda\nabla_j \varphi_{\mu k \lambda} \varphi_{iab} \o^{ka} \o^{pb}
&=&
2|\varphi|^2  N^\mu{}_p{}^\lambda (- {1 \over 2} \alpha_j \varphi_{\mu k\lambda} + {1 \over 2} \alpha_{J j} \varphi_{J \mu, k\lambda} - E_{j;\mu k \lambda}) \varphi_{iab} \o^{ka} \o^{pb}.
\nonumber
\eea
We work out the bilinears term by term
\bea
2|\varphi|^2  N^\mu{}_p{}^\lambda (- {1 \over 2} \alpha_j \varphi_{\mu k\lambda})\varphi_{iab} \o^{ka} \o^{pb}&=& 
-|\varphi|^2  \alpha_j  N^\mu{}_p{}^\lambda \varphi_{\mu k\lambda}\varphi_{iab} \o^{ka} \o^{pb}\nonumber\\
&=&
-{1\over 4}|\varphi|^4  \alpha_j  N^\mu{}_p{}^\lambda (\o_{\mu i} g_{\lambda b} + \o_{\lambda b} g_{\mu i} - \o_{\mu b} g_{\lambda i} - \o_{\lambda i} g_{\mu b}) \o^{pb}\nonumber\\
&=&
-{1\over 4}|\varphi|^4  \alpha_j  N^\mu{}_p{}^\lambda (\o_{\mu i} J^p{}_\lambda - \delta^p{}_\lambda g_{\mu i} + \delta^p{}_\mu g_{\lambda i} - \o_{\lambda i} J^p{}_\mu)\nonumber\\
&=&
-{1\over 4}|\varphi|^4  \alpha_j (- N_{Ji, p}{}^{Jp} - N_{ip}{}^p + N^p{}_{pi} + N^p{}_{p, Ji})=0. \nonumber
\eea
\bea
2|\varphi|^2  N^\mu{}_p{}^\lambda  ({1\over 2} \alpha_{J j} \varphi_{J \mu, k\lambda}) \varphi_{iab} \o^{ka} \o^{pb}&=& |\varphi|^2 \alpha_{J j} N^\mu{}_p{}^\lambda   \varphi_{J \mu, k\lambda} \varphi_{iab} \o^{ka} \o^{pb}\nonumber\\
&=&
-{1\over 4}|\varphi|^4  \alpha_{Jj}  N^\mu{}_p{}^\lambda (\o_{J\mu, i} g_{\lambda b} + \o_{\lambda b} g_{J\mu, i} - \o_{J\mu, b} g_{\lambda i} - \o_{\lambda i} g_{J\mu, b}) \o^{pb}\nonumber\\
&=&
-{1\over 4}|\varphi|^4  \alpha_{Jj}  N^\mu{}_p{}^\lambda (- g_{\mu i} J^p{}_\lambda- \delta^p{}_\lambda \o_{\mu i } + J^p{}_{\mu} g_{\lambda i} + \delta^p{}_\mu \o_{\lambda i} )\nonumber\\
&=&
-{1\over 4}|\varphi|^4  \alpha_{Jj}(- N_{i p}{}^{Jp} + N_{Ji, p}{}^p + N^{Jp}{}_{p i} - N^p{}_{p, Ji})=0.\nonumber
\eea
\bea\label{E-terms}
2|\varphi|^2  N^\mu{}_p{}^\lambda (- E_{j;\mu k \lambda}) \varphi_{iab} \o^{ka} \o^{pb}&=& 
-2|\varphi|^2  N^\mu{}_p{}^\lambda E_{j;\mu k \lambda}\varphi_{iab} \o^{ka} \o^{pb}\nonumber\\
&=& 
-2|\varphi|^2  N^\mu{}_p{}^\lambda( N_{j\mu}{}^\ell \varphi_{\ell k \lambda}+ N_{jk}{}^\ell \varphi_{\mu \ell \lambda}+ N_{j\lambda}{}^\ell \varphi_{\mu k \ell })\varphi_{iab} \o^{ka} \o^{pb}\nonumber\\
\eea
The first term in the above last line is easy to handle
\bea
2|\varphi|^2  N^\mu{}_p{}^\lambda N_{j\mu}{}^\ell \varphi_{\ell k \lambda}\varphi_{iab} \o^{ka} \o^{pb}
&=&
{1\over 2}|\varphi|^4  N^\mu{}_p{}^\lambda N_{j\mu}{}^\ell (\o_{\ell i } g_{\lambda b} + \o_{\lambda b} g_{\ell i} - \o_{\ell b} g_{\lambda i} - \o_{\lambda i} g_{\ell b} ) \o^{pb}\nonumber\\
&=&
{1\over 2}|\varphi|^4  N^\mu{}_p{}^\lambda N_{j\mu}{}^\ell (\o_{\ell i} J^p{}_\lambda - \delta^p{}_\lambda g_{\ell i} + \delta^p{}_\ell g_{\lambda i} - \o_{\lambda i} J^p{}_{\ell})\nonumber\\
&=&{1\over 2}|\varphi|^4 (-N^\mu{}_{J\lambda}{}^\lambda N_{j \mu, Ji}- N^\mu{}_{p}{}^p N_{j\mu i} + N^\mu{}_{pi} N_{j\mu}{}^p + N^\mu{}_{J\ell, Ji} N_{j \mu}{}^\ell\nonumber\\
&=&{1\over 2}|\varphi|^4(N^\mu{}_{pi} N_{j\mu}{}^p - N^\mu{}_{pi} N_{j \mu}{}^p)=0.
\eea
The third term can also be handled in the similar way
\bea
2|\varphi|^2  N^\mu{}_p{}^\lambda N_{j\lambda}{}^\ell \varphi_{\mu k \ell }\varphi_{iab} \o^{ka} \o^{pb}
&=& {1\over 2}|\varphi|^4  N^\mu{}_p{}^\lambda N_{j\lambda}{}^\ell ( \o_{\mu i}g_{\ell b} + \o_{\ell b} g_{\mu i} - \o_{\mu b} g_{\ell i} - \o_{\ell i} g_{\mu b} ) \o^{pb}\nonumber\\
&=&
{1\over 2}|\varphi|^4  N^\mu{}_p{}^\lambda N_{j\lambda}{}^\ell ( \o_{\mu i} J^p{}_\ell - \delta^p{}_\ell g_{\mu i} + \delta^p{}_\mu g_{\ell i} - \o_{\ell i} J^p{}_\mu)\nonumber\\
&=&
{1\over 2}|\varphi|^4 (- N_{Ji, p}{}^\lambda N_{j\lambda}{}^{Jp}- N_{ip}{}^\lambda N_{j\lambda}{}^p + N^p{}_p{}^\lambda N_{j\lambda i} + N^{Jp}{}_p{}^\lambda N_{j\lambda, Ji})\nonumber\\
&=&
{1\over 2}|\varphi|^4 (N_{i p}{}^\lambda N_{j\lambda}{}^{p}- N_{ip}{}^\lambda N_{j\lambda}{}^p)=0.
\eea
For the second term in (\ref{E-terms}), we will use the Bianchi identity and switch the indices as before, 
$N^p{}_{ij} \varphi_{pkl} = - N^p{}_{kl} \varphi_{p ij}$, obtaining
\bea
N^\mu{}_p{}^\lambda N_{jk}{}^\ell \varphi_{\mu \ell \lambda}&=&
- N^\mu{}_p{}^\lambda(N^\ell{}_{jk} + N_{k}{}^\ell{}_j) \varphi_{\mu \ell \lambda} 
=
N^\mu{}_p{}^\lambda\, N^\ell{}_{jk}\varphi_{\ell \mu \lambda} - N_{k}{}^\ell{}_j\, N^\mu{}_p{}^\lambda \varphi_{\mu \ell \lambda}\nonumber\\
&=& 
- N^\mu{}_p{}^\lambda\, N^\ell{}_{\mu \lambda} \varphi_{\ell j k} + N_{k}{}^\ell{}_j\, N^\mu{}_{\ell}{}^\lambda \varphi_{\mu p \lambda}
\eea
Therefore, 
\bea
-2|\varphi|^2  N^\mu{}_p{}^\lambda N_{jk}{}^\ell \varphi_{\mu \ell \lambda}\varphi_{iab} \o^{ka} \o^{pb}&=&
2|\varphi|^2 ( N^\mu{}_p{}^\lambda\, N^\ell{}_{\mu \lambda} \varphi_{\ell j k} - N_{k}{}^\ell{}_j\, N^\mu{}_{\ell}{}^\lambda \varphi_{\mu p \lambda}) \varphi_{iab} \o^{ka} \o^{pb}\nonumber\\
&=&
-{1\over 2}|\varphi|^4 N^\mu{}_p{}^\lambda N^\ell{}_{\mu \lambda} ( \o_{\ell i} g_{jb} + \o_{j b} g_{\ell i} - \o_{\ell b} g_{j i } - \o_{j i} g_{\ell b} )\o^{pb}  \nonumber\\
&&
+ {1\over 2}|\varphi|^4 N_{k}{}^\ell{}_j\, N^\mu{}_{\ell}{}^\lambda( \o_{\mu i} g_{\lambda a} + \o_{\lambda a } \o_{\mu i} - \o_{\mu a} g_{\lambda i} - \o_{\lambda i}g_{\mu a} ) \o^{ka}
\nonumber\\
&=&
-{1\over 2}|\varphi|^4 N^\mu{}_p{}^\lambda N^\ell{}_{\mu \lambda} ( \o_{\ell i} J^p{}_j - \delta^p{}_j g_{\ell i} + \delta^p{}_\ell g_{j i } - \o_{j i} J^p{}_\ell )  \nonumber\\
&&
\quad+{1\over 2}|\varphi|^4N_{k}{}^\ell{}_j\, N^\mu{}_{\ell}{}^\lambda( \o_{\mu i} J^k{}_\lambda - \delta^k{}_\lambda g_{\mu i} + \delta^k{}_\mu g_{\lambda i} - \o_{\lambda i}J^k{}_\mu ) \nonumber\\
&=&
-{1\over 2}|\varphi|^4(- N^\mu{}_{Jj}{}^\lambda N_{Ji, \mu \lambda}- N^\mu{}_j{}^\lambda N_{i \mu \lambda} + N^\mu{}_p{}^\lambda N^p{}_{\mu \lambda}g_{ji} )\nonumber\\
&&
+{1\over 2}|\varphi|^4( -N_{J\lambda}{}^\ell{}_j\, N_{Ji, \ell}{}^\lambda-N_{\lambda}{}^\ell{}_j\, N_{i \ell}{}^\lambda + N_{\mu}{}^\ell{}_j\, N^\mu{}_{\ell i}
+ N_{J\mu}{}^\ell{}_j\, N^\mu{}_{\ell, Ji})\nonumber
\eea
The right hand side can be readily simplified as follows,
\bea
&&
-{1\over 2}|\varphi|^4 ( - 2N^\mu{}_j{}^\lambda N_{i \mu \lambda} + N^\mu{}_p{}^\lambda N^p{}_{\mu \lambda}g_{ji} )
+{1\over 2}|\varphi|^4(-2N_{\lambda}{}^\ell{}_j\, N_{i\ell}{}^\lambda + N_{\mu}{}^\ell{}_j\, N^\mu{}_{\ell i}
- N_{J\mu}{}^\ell{}_j\, N^{J\mu}{}_{\ell i})\nonumber\\
&=& 
-{1\over 2}|\varphi|^4 ( - 2N^\mu{}_j{}^\lambda N_{i \mu \lambda} + N^\mu{}_p{}^\lambda N^p{}_{\mu \lambda}g_{ji}+ 2N_{\lambda}{}^\ell{}_j\, N_{i\ell}{}^\lambda -2 N_{\mu}{}^\ell{}_j\, N^\mu{}_{\ell i})\nonumber\\
&=&
-{1\over 2}|\varphi|^4(2N^{\mu\lambda}{}_j N_{i \mu \lambda} - 2N^{\lambda\ell}{}_j\, N_{i\lambda \ell}+  (N_{-}^2)^\lambda{}_{\lambda}g_{ji}-2 N_{\mu}{}^\ell{}_j\, N^\mu{}_{\ell i})\nonumber\\
&=& 
-|\varphi|^4 ({1\over 2} (N_{-}^2)^\lambda{}_{\lambda} g_{ij} - N_{\mu}{}^\ell{}_j\, N^\mu{}_{\ell i}).
\eea
Putting the above computations into (\ref{3rd-term}), we obtain
\bea
2|\varphi|^2  N^\mu{}_p{}^\lambda\nabla_j \varphi_{\mu k \lambda} \varphi_{iab} \o^{ka} \o^{pb}=-|\varphi|^4 ({1\over 2}  (N_{-}^2)^\lambda{}_{\lambda} g_{ij} - N_{\mu}{}^\ell{}_j\, N^\mu{}_{\ell i})
\eea
Therefore, we obtain the first term (I) in (\ref{Ndagger}):
\bea
{\rm (I)}\cdot \varphi_{iab}\o^{ka}\o^{pb}&=& |\varphi|^4(N_{ip}{}^\lambda N_{j\lambda}{}^p-N^\lambda{}_{pi} \, N_{j \lambda}{}^{p})- |\varphi|^4 ({1\over 2}  (N_{-}^2)^\lambda{}_{\lambda} g_{ij} - N_{p}{}^\lambda{}_j\, N^p{}_{\lambda i})\nonumber\\
&=&
-{1\over 2} |\varphi|^4  (N_{-}^2)^\lambda{}_{\lambda} g_{ij}+ |\varphi|^4(N_{ip}{}^\lambda N_{j\lambda}{}^p-N^\lambda{}_{pi} \, N_{j \lambda}{}^{p}+N_{p}{}^\lambda{}_j\, N^p{}_{\lambda i})\nonumber
\eea
Using the Bianchi identity, we readily find
\bea
N_{ip}{}^\lambda N_{j\lambda}{}^p-N^\lambda{}_{pi} \, N_{j \lambda}{}^{p}+N_{p}{}^\lambda{}_j\, N^p{}_{\lambda i}=- N_{\lambda j}{}^p N_p{}^\lambda{}_i
\eea
%
Therefore, 
\bea
{\rm (I)}\cdot \varphi_{iab}\o^{ka}\o^{pb}&=&-{1\over 2} |\varphi|^4  (N_{-}^2)^\lambda{}_{\lambda} g_{ij}+ |\varphi|^4 N_{\lambda p j} N^{p\lambda}{}_i.
\eea


\subsubsection{Computation for (II)}



Next we work out the contributions of (II) in (\ref{Ndagger}). The contributions from (III) will turn out to be similar.
\bea\label{Ndagger2}
{1\over 2} {\rm II}&=&{1\over 2} \nabla_p (|\varphi|^2(N^\dagger\cdot\varphi)_{jk})=-{1\over 2} \nabla_p (|\varphi|^2(N^\dagger\cdot\varphi)_{kj}) \nonumber\\
&=& 
 |\varphi|^2 \alpha_p N^\mu{}_j{}^\lambda \varphi_{\mu k \lambda}  - |\varphi|^2\nabla_p N^\mu{}_j{}^\lambda \varphi_{\mu k \lambda}-|\varphi|^2 N^\mu{}_j{}^\lambda\nabla_p \varphi_{\mu k \lambda}
\eea
Again, we will work out the bilinears term by term.
\bea\label{term-1}
|\varphi|^2 \alpha_p N^\mu{}_j{}^\lambda \varphi_{\mu k \lambda} \varphi_{iab}\o^{ka}\o^{pb}
&=&
 {1\over 4} |\varphi|^4 \alpha_p N^\mu{}_j{}^\lambda (\o_{\mu i} g_{\lambda b} + \o_{\lambda b} g_{\mu i} - \o_{\mu b} g_{\lambda i} - \o_{\lambda i} g_{\mu b}) \o^{pb}\nonumber\\
 &=&
 {1\over 4}   |\varphi|^4 \alpha_p N^\mu{}_j{}^\lambda (\o_{\mu i} J^p{}_\lambda - \delta^p{}_\lambda g_{\mu i} + \delta^p{}_\mu g_{\lambda i} - \o_{\lambda i} J^p{}_\mu)\nonumber\\
 &=&
 {1\over 4}  |\varphi|^4 \alpha_p( -N_{Ji, j}{}^{Jp} - N_{ij}{}^p + N^p{}_{ji} + N^{Jp}{}_{j, Ji})\nonumber\\
 &=&
 {1\over 4}  |\varphi|^4 \alpha_p(- N_{ij}{}^p- N_{ij}{}^p+N^p{}_{ji} +N^p{}_{ji} )\nonumber\\
 &=&
 {1\over 2} |\varphi|^4 \alpha_p(- N_{ij}{}^p +N^p{}_{ji} )
\eea
Next, we deal with the second term in (\ref{Ndagger2})
\bea\label{2nd-J}
-|\varphi|^2\nabla_p N^\mu{}_j{}^\lambda \varphi_{\mu k\lambda} \varphi_{iab}\o^{ka}\o^{pb}&=&
 -{1\over 4} |\varphi|^4\nabla_p N^\mu{}_j{}^\lambda (\o_{\mu i} g_{\lambda b} + \o_{\lambda b} g_{\mu i} - \o_{\mu b} g_{\lambda i} - \o_{\lambda i} g_{\mu b}) \o^{pb}\nonumber\\
 &=&
 -{1\over 4} |\varphi|^4\nabla_p N^\mu{}_j{}^\lambda (\o_{\mu i} J^p{}_\lambda - \delta^p{}_\lambda g_{\mu i} + \delta^p{}_\mu g_{\lambda i} - \o_{\lambda i} J^p{}_\mu)\nonumber\\
 &=&
 {1\over 4} |\varphi|^4( \nabla_p  N_{ij}{}^p - \nabla_p N^p{}_{ji} )\nonumber\\
 &&- {1\over 4} |\varphi|^4 (\o_{\mu i}\nabla_p N^\mu{}_j{}^\lambda J^p{}_\lambda - \o_{\lambda i} \nabla_p N^\mu{}_j{}^\lambda J^p{}_\mu)
\eea
For the second group of terms in (\ref{2nd-J}) , we need to take care of $\nabla J$,
\bea
&&\o_{\mu i}\nabla_p N^\mu{}_j{}^\lambda J^p{}_\lambda - \o_{\lambda i} \nabla_p N^\mu{}_j{}^\lambda J^p{}_\mu\\ 
&=& 
\o_{\mu i}\nabla_p N^\mu{}_j{}^\lambda J^p{}_\lambda - \o_{\mu i} \nabla_p N^\lambda{}_j{}^\mu J^p{}_\lambda 
=
\o_{\mu i}\nabla_p( N^\mu{}_j{}^\lambda + N^{\lambda\mu}{}_j) J^p{}_\lambda 
=
-\o_{\mu i}\nabla_p N_j{}^{\lambda \mu}J^p{}_\lambda\nonumber\\
&=&
-\o_{\mu i}(\nabla_p( N_j{}^{\lambda \mu}J^p{}_\lambda)  - N_j{}^{\lambda \mu}\nabla_p J^p{}_\lambda)
=
 \o_{\mu i} \nabla_p N_{Jj}{}^{p\mu} - 2 \o_{\mu i}N_j{}^{\lambda \mu}N_{p\lambda}{}^{Jp}\nonumber\\
 &=& \o_{\mu i} \nabla_p N_{Jj}{}^{p\mu} 
 = 
 \nabla_p N_{Jj}{}^{p\mu} J^\ell{}_\mu g_{\ell i} 
 =
  \nabla_p(N_{Jj}{}^{p\mu} J^\ell{}_\mu g_{\ell i}) - N_{Jj}{}^{p\mu} \nabla_p J^\ell{}_\mu g_{\ell i}\nonumber\\
 &=&
 \nabla_p N_j{}^p{}_i - 2N_j{}^{p\mu} N_{p\mu i}
\eea
Putting this back into the calculation,
\bea
-|\varphi|^2\nabla_p N^\mu{}_j{}^\lambda \varphi_{\mu k\lambda} \varphi_{iab}\o^{ka}\o^{pb}
&=&
{1\over 4} |\varphi|^4( \nabla_p  N_{ij}{}^p - \nabla_p N^p{}_{ji} - \nabla_p N_j{}^p{}_i )
+{1\over 2} |\varphi|^4 N_j{}^{p\mu} N_{p\mu i}.\nonumber\\
&=&
{1\over 2} |\varphi|^4( \nabla_p  N_{ij}{}^p - \nabla_p N^p{}_{ji}  )
+ {1\over 2} |\varphi|^4 N_j{}^{p\mu} N_{p\mu i}
\eea
where we used the Bianchi identity $- N_j{}^p{}_i = N_{ij}{}^p +N^p{}_{ij} $ to obtain the last equality above.

Now, we deal with the $\nabla N$ terms using the projected Levi-Civita connection
\bea
\nabla_p  N_{ij}{}^p - \nabla_p N^p{}_{ji} &=& 
{\frak D}_p N_{ij}{}^p - N_{\alpha j}{}^p N_{pi}{}^\alpha - N_{i\alpha}{}^p N_{pj}{}^\alpha - N_{ij\alpha} N_p{}^{p\alpha}\\
&&
- {\frak D}_p N^p{}_{ji} + N_{\alpha ji} N_p{}^{p\alpha} + N^p{}_{\alpha i} N_{pj}{}^\alpha + N^p{}_{j\alpha} N_{pi}{}^\alpha\nonumber\\
&=& 
{\frak D}_p N_{ij}{}^p- {\frak D}_p N^p{}_{ji}- (N_{\alpha j}{}^p - N^p{}_{j\alpha}) N_{pi}{}^\alpha -(N_{i\alpha}{}^p-N^p{}_{\alpha i}) N_{pj}{}^\alpha \nonumber
\eea
since $N_p{}^{p\alpha}=0$. Next, apply the Bianchi identity of $N$ to the last two terms, and get
\bea
\nabla_p  N_{ij}{}^p - \nabla_p N^p{}_{ji} = {\frak D}_p N_{ij}{}^p- {\frak D}_p N^p{}_{ji} + N_j{}^p{}_\alpha N_{pi}{}^\alpha + N_\alpha{}^p{}_i N_{pj}{}^\alpha.
\eea
So, we have 
\bea\label{term-2}
-|\varphi|^2\nabla_p N^\mu{}_j{}^\lambda \varphi_{\mu k\lambda} \varphi_{iab}\o^{ka}\o^{pb}
&=& 
 {1\over 2} |\varphi|^4 ({\frak D}_p N_{ij}{}^p- {\frak D}_p N^p{}_{ji}) \\
&&
+ {1\over 2} |\varphi|^4  (N_j{}^p{}_\alpha N_{pi}{}^\alpha + N_\alpha{}^p{}_i N_{pj}{}^\alpha+  N_j{}^{p\alpha} N_{p\alpha i})\nonumber\\
&=&
{1\over 2} |\varphi|^4 ({\frak D}_p N_{ij}{}^p- {\frak D}_p N^p{}_{ji}) +{1\over 2} |\varphi|^4 N_\alpha{}^p{}_i N_{pj}{}^\alpha \nonumber
\eea
Next, we deal with the last term in (\ref{Ndagger2}).
Since
$\nabla_p \varphi_{\mu k\lambda} = - {1 \over 2} \alpha_p \varphi_{\mu k\lambda} + {1 \over 2} \alpha_{J p} \varphi_{J \mu, k\lambda} - E_{p;\mu k \lambda}$,
we have
\bea
-|\varphi|^2 N^\mu{}_j{}^\lambda\nabla_p \varphi_{\mu k \lambda}&=&
|\varphi|^2 N^\mu{}_j{}^\lambda({1 \over 2} \alpha_p \varphi_{\mu k\lambda} - {1 \over 2} \alpha_{J p} \varphi_{J \mu, k\lambda} + E_{p;\mu k \lambda}).
\eea
We work out the bilinears term by term.
\bea
{1 \over 2}|\varphi|^2 \alpha_p N^\mu{}_j{}^\lambda \varphi_{\mu k\lambda}\varphi_{iab} \o^{ka} \o^{pb}
&=&
{1 \over 8}|\varphi|^4 \alpha_p N^\mu{}_j{}^\lambda (\o_{\mu i} J^p{}_\lambda - \delta^p{}_\lambda g_{\mu i} + \delta^p{}_\mu g_{\lambda i} - \o_{\lambda i} J^p{}_\mu)\nonumber\\
 &=&
 {1 \over 8}|\varphi|^4 \alpha_p (- N_{Ji, j}{}^{Jp} - N_{ij}{}^p + N^p{}_{ji} + N^{Jp}{}_{j, Ji})\nonumber\\
 &=&
 {1 \over 4}|\varphi|^4 \alpha_p (- N_{ij}{}^p + N^p{}_{ji})
\eea
\bea
-{1 \over 2}|\varphi|^2 \alpha_{Jp} N^\mu{}_j{}^\lambda \varphi_{J\mu, k\lambda}\varphi_{iab} \o^{ka} \o^{pb}
&=&
-{1 \over 8}|\varphi|^4 \alpha_{Jp} N^\mu{}_j{}^\lambda (- g_{\mu i} J^p{}_\lambda- \delta^p{}_\lambda \o_{\mu i } + J^p{}_{\mu} g_{\lambda i} + \delta^p{}_\mu \o_{\lambda i} )\nonumber\\
&=&
-{1 \over 8}|\varphi|^4 \alpha_{Jp} (-N_{ij}{}^{Jp} + N_{Ji, j}{}^p + N^{Jp}{}_{j i} - N^p{}_{j, Ji})\nonumber\\
&=&
-{1 \over 8}|\varphi|^4 \alpha_{Jp} (-2N_{ij}{}^{Jp}+ 2N^{Jp}{}_{j i})
={1 \over 4}|\varphi|^4 \alpha_{p} N_{j}{}^p{}_i
\eea
by the Bianchi identity satisfied by $N$. The terms $E$ lead to
\bea\label{Eterm}
|\varphi|^2 N^\mu{}_j{}^\lambda E_{p;\mu k \lambda} \varphi_{iab} \o^{ka} \o^{pb}
= |\varphi|^2 N^\mu{}_j{}^\lambda ( N_{p\mu}{}^\ell \varphi_{\ell k \lambda}+ N_{pk}{}^\ell \varphi_{\mu \ell \lambda}+ N_{p\lambda}{}^\ell \varphi_{\mu k \ell })\varphi_{iab} \o^{ka} \o^{pb}\nonumber
\eea
We compute the three terms
\bea
|\varphi|^2 N^\mu{}_j{}^\lambda N_{p\mu}{}^\ell \varphi_{\ell k \lambda}\varphi_{iab} \o^{ka} \o^{pb}
&=&
{1\over 2}|\varphi|^4 N^\mu{}_j{}^\lambda N_{p\mu}{}^\ell ( \o_{\ell i} g_{\lambda b} + \o_{\lambda  b} g_{\ell i} - \o_{\ell b} g_{\lambda i} - \o_{\lambda i } g_{\ell b} ) \o^{pb}\nonumber\\
&=&
{1\over 2}|\varphi|^4 N^\mu{}_j{}^\lambda N_{p\mu}{}^\ell (\o_{\ell i} J^p{}_\lambda - \delta^p{}_\lambda g_{\ell i} + \delta^p{}_\ell g_{\lambda i} - \o_{\lambda i} J^p{}_\ell)\nonumber\\
&=&
{1\over 2}|\varphi|^4(-N^\mu{}_j{}^{Jp} N_{p\mu, Ji} - N^\mu{}_j{}^p N_{p\mu i}+ N^\mu{}_{ji} N_{p\mu}{}^p + N^\mu{}_{j, Ji} N_{p\mu}{}^{Jp} )\nonumber\\
&=&
{1\over 2}|\varphi|^4(N^\mu{}_j{}^{p} N_{p\mu i}-N^\mu{}_j{}^p N_{p\mu i})=0
\eea
\bea
 |\varphi|^2 N^\mu{}_j{}^\lambda N_{p\lambda}{}^\ell \varphi_{\mu k \ell }\varphi_{iab} \o^{ka} \o^{pb}
 &=&
 {1\over 2} |\varphi|^4 N^\mu{}_j{}^\lambda N_{p\lambda}{}^\ell
  (\o_{\mu i} g_{\ell b} + \o_{\ell b} g_{\mu i } - \o_{\mu b} g_{\ell i} - \o_{\ell i} g_{\mu b} ) \o^{pb} \nonumber\\
&=&
{1\over 2}|\varphi|^4 N^\mu{}_j{}^\lambda N_{p\lambda}{}^\ell( \o_{\mu i } J^p{}_\ell - \delta^p{}_\ell g_{\mu i} + \delta^p{}_\mu g_{\ell i} - \o_{\ell i} J^p{}_\mu)  \nonumber\\
&=&
{1\over 2}|\varphi|^4 (-N_{Ji, j}{}^\lambda N_{p\lambda}{}^{Jp}-N_{ij}{}^\lambda N_{p\lambda}{}^p+ N^p{}_j{}^\lambda N_{p\lambda i}+ N^{Jp}{}_j{}^\lambda N_{p\lambda, Ji}) \nonumber\\
&=&
{1\over 2}|\varphi|^4(N^p{}_j{}^\lambda N_{p\lambda i}- N^{Jp}{}_j{}^\lambda N_{p\lambda, Ji})=0
\eea
The second term in (\ref{Eterm}) is more complicated, we first note that, by interchanging indices $k\leftrightarrow p$ and $a\leftrightarrow b$, 
\bea
N^\mu{}_j{}^\lambda N_{pk}{}^\ell \varphi_{\mu \ell \lambda}\varphi_{iab} \o^{ka} \o^{pb}
&=& -N^\mu{}_j{}^\lambda N_{kp}{}^\ell \varphi_{\mu \ell \lambda}\varphi_{iab} \o^{ka} \o^{pb}
\eea
It follows that
\bea
N^\mu{}_j{}^\lambda N_{pk}{}^\ell \varphi_{\mu \ell \lambda}\varphi_{iab} \o^{ka} \o^{pb}
&=&
 {1\over 2}N^\mu{}_j{}^\lambda (N_{pk}{}^\ell -N_{kp}{}^\ell )\varphi_{\mu \ell \lambda}\varphi_{iab} \o^{ka} \o^{pb} \nonumber\\
 &=& 
 {1\over 2}N^\mu{}_j{}^\lambda (N_{pk}{}^\ell +N_{k}{}^\ell{}_p )\varphi_{\mu \ell \lambda}\varphi_{iab} \o^{ka} \o^{pb} \nonumber\\
 &=& -{1\over 2} N^\mu{}_j{}^\lambda N^\ell{}_{pk} \varphi_{\mu \ell \lambda}\varphi_{iab} \o^{ka} \o^{pb}
\eea
where we use Bianchi identity to get the last equality. Now, we can ready to use the identity $N^{p}{}_{k\ell} \varphi_{pij} = - N^p{}_{ij} \varphi_{pk\ell}$ to handle the second term in (\ref{Eterm})
\bea
|\varphi|^2 N^\mu{}_j{}^\lambda N_{pk}{}^\ell \varphi_{\mu \ell \lambda}\varphi_{iab} \o^{ka} \o^{pb}
&=&
-{1\over 2} |\varphi|^2N^\mu{}_j{}^\lambda N^\ell{}_{pk} \varphi_{\mu \ell \lambda}\varphi_{iab} \o^{ka} \o^{pb} \nonumber\\
 &=& 
 {1\over 2} |\varphi|^2N^\mu{}_j{}^\lambda N^\ell{}_{\mu \lambda} \varphi_{\ell kp} \varphi_{iab} \o^{ka} \o^{pb} \nonumber\\
 &=& 
 -{1\over 2} |\varphi|^4 N^\mu{}_j{}^\lambda N^\ell{}_{\mu \lambda} g_{\ell i} 
 =
  {1\over 2} |\varphi|^4 N^{\mu\lambda}{}_j N_{i \mu \lambda}.
 \eea
So,
\bea
|\varphi|^2 N^\mu{}_j{}^\lambda E_{p;\mu k \lambda} \varphi_{iab} \o^{ka} \o^{pb}= 
 {1\over 2} |\varphi|^4 N^{\mu\lambda}{}_j N_{i \mu \lambda}.
 \eea
Putting the above calculation together, we obtain
\bea\label{term-3}
-|\varphi|^2 N^\mu{}_j{}^\lambda\nabla_p \varphi_{\mu k \lambda}\varphi_{iab} \o^{ka} \o^{pb}&=& {1 \over 4}|\varphi|^4 \alpha_p (- N_{ij}{}^p + N^p{}_{ji})+{1 \over 4}|\varphi|^4 \alpha_{p} N_{j}{}^p{}_i+  {1\over 2} |\varphi|^4 N^{\mu\lambda}{}_j N_{i \mu \lambda} \nonumber\\
& =&  {1\over 2} |\varphi|^4 N^{\mu\lambda}{}_j N_{i \mu \lambda} +{1 \over 2}|\varphi|^4 \alpha_{p} N_{j}{}^p{}_i
\eea
using Bianchi identity $N^p{}_{ji} + N_i{}^p{}_j + N_{ji}{}^p=0$.

Back to (\ref{Ndagger2}), using (\ref{term-1}) (\ref{term-2}) and (\ref{term-3}), we complete the calculation for (II):
\bea
{\rm (II)}\cdot \varphi_{iab}\o^{ka}\o^{pb}&=& |\varphi|^4 \alpha_p(- N_{ij}{}^p +N^p{}_{ji} )
+  |\varphi|^4 ({\frak D}_p N_{ij}{}^p- {\frak D}_p N^p{}_{ji})\\
&&
+  |\varphi|^4  N_\lambda{}^p{}_i N_{pj}{}^\lambda+ |\varphi|^4 N^{\mu\lambda}{}_j N_{i \mu \lambda}+|\varphi|^4 \alpha_{p} N_{j}{}^p{}_i
\nonumber\\
&=&
  |\varphi|^4 {\frak D}_p N_{ji}{}^p+ 2 |\varphi|^4 \alpha_p N_j{}^p{}_i
- |\varphi|^4 N^{p\lambda}{}_i N_{p\lambda j}
\nonumber
\eea
Note that $N^p{}_{ji}=0$ up to the symmetrization for $(i\leftrightarrow j)$. So, terms involving $N^p{}_{ij}$ vanish up to the symmetrization for $ (i\leftrightarrow j)$. For the two quadratic terms about $N$, we use Bianchi identity to obtain the last line.
Thus
\bea
 {\rm (II)}\cdot \varphi_{iab}\o^{ka}\o^{pb}= |\varphi|^4\left\{ {\frak D}_p N_{ji}{}^p+ 2\alpha_p N_j{}^p{}_i  -  N^{p\lambda}{}_i N_{p\lambda j}\right\}
\eea

\subsubsection{Computation for (III)}
Next, we consider (III) in (\ref{Ndagger}). We simply observe that by switching the indices $k\leftrightarrow p$ and $a\leftrightarrow b$ and exploiting the antisymmetry of $(N^\dagger\varphi)_{kj}$ in $j$ and $k$, we may write
\bea
{\rm (III)}\cdot \varphi_{iab}\o^{ka}\o^{pb} &=& \nabla_k (|\varphi|^2(N^\dagger\cdot\varphi)_{pj}) \varphi_{iab}\o^{ka}\o^{pb} \nonumber\\
&=& - \nabla_p (|\varphi|^2(N^\dagger\cdot\varphi)_{kj})\varphi_{iab}\o^{ka}\o^{pb}\nonumber\\
&=& \nabla_p (|\varphi|^2(N^\dagger\cdot\varphi)_{jk})\varphi_{iab}\o^{ka}\o^{pb} =
{\rm (II)}\cdot \varphi_{iab}\o^{ka}\o^{pb}.
\eea
We can now put (I), (II) and (III) all together, 
\bea
(d(|\varphi|^2N^\dagger\cdot\varphi))_{jkp}\cdot \varphi_{iab}\o^{ka}\o^{pb}+ (i\leftrightarrow j)
&=&
- |\varphi|^4  (N_{-}^2)^\lambda{}_{\lambda} g_{ij}+ |\varphi|^4\left\{N_{\lambda p j} N^{p\lambda}{}_i+(i\leftrightarrow j)\right\}\nonumber\\
&& + 2|\varphi|^4\big\{ {\frak D}_p N_{ji}{}^p+ 2\alpha_p N_j{}^p{}_i  -  N^{p\lambda}{}_i N_{p\lambda j}+ (i\leftrightarrow j)\big\}\nonumber
 \eea

%
\begin{lemma} \label{dNdagger} In conclusion, we have
\bea
&&d(|\varphi|^2N^\dagger\cdot\varphi)_{jkp}\cdot \varphi_{iab}\o^{ka}\o^{pb}+ (i\leftrightarrow j)
\\
&=&
|\varphi|^4\big\{ 2({\frak D}_p N_{ji}{}^p + {\frak D}_p N_{ij}{}^p) + 4 \alpha_p (N_j{}^p{}_i + N_i{}^p{}_j) -  (N_{-}^2)^\lambda{}_{\lambda} g_{ij} - 4(N_+^2)_{ij} + 2(N_{-}^2)_{ij} \big\}\nonumber
\eea
\end{lemma}

\subsection{The flow of $g_{ij}$}
Assembling all the terms in (\ref{phi-flow-rhs}) and putting them in (\ref{flow-tildeg}), we obtain the flow of $\tilde g_{ij}$,
\bea
\partial_t \tilde{g}_{ij} &=& - \bigg\{ (-| \varphi |^2 d d^\dagger \varphi)_{jkp} \varphi_{iab}  \omega^{ka} \omega^{pb} - (d | \varphi |^2 \wedge d^\dagger \varphi)_{jkp} \varphi_{iab}  \omega^{ka} \omega^{pb}  \\
&&+ (d \iota_{ \nabla | \varphi |^2} \varphi)_{jkp} \varphi_{iab}  \omega^{ka} \omega^{pb} + (2d(| \varphi |^2 N^\dagger \cdot \varphi))_{jkp} \varphi_{iab}  \omega^{ka} \omega^{pb} + (i \leftrightarrow j) \bigg\} \nonumber
\eea
By (\ref{bochner-contribution}), (\ref{grad-dagger}), Lemmas \ref{d-iota-contribution} and \ref{dNdagger}, and the identity (\ref{qua-N}),
\bea
\partial_t \tilde{g}_{ij} &=& - | \varphi|^4 \bigg\{ 2 ({\frak D}_k N_{ij}{}^k  +  {\frak D}_k N_{ji}{}^k) + R g_{ij}
+  2 \nabla_\mu \alpha^\mu g_{ij}  + {1 \over 2}( \nabla_j \alpha_i + \nabla_i \alpha_j)
\\
&& - {1 \over 2} (J^p{}_j J^q{}_i \nabla_{p} \alpha_{q} + J^p{}_i J^q{}_j \nabla_{p} \alpha_{q}) 
- \alpha_i \alpha_j  + \alpha_{Ji} \alpha_{Jj}  -  2 \alpha_\mu \alpha^\mu g_{ij}  +4 \alpha_p (N_j{}^p{}_i + N_i{}^p{}_j)\bigg\}
\nonumber
\eea
Recall that $\tilde{g}_{ij} = | \varphi |^2 g_{ij}$. Therefore
\bea
\p_t \log \det \tilde{g} = | \varphi |^{-2} g^{ij} \partial_t \tilde{g}_{ij} =
| \varphi|^2 \big\{ -12 \nabla_\mu \alpha^\mu - 6 R + 12 |\alpha|^2  \big\} .
\eea
Since ${\rm det} \tilde{g} = |\varphi|^{12} {\rm det} g$ and $\partial_t {\rm det} \, g = 0$ as the volume form of $g$ equals to $\o^3/3!$ and $\omega$ is fixed, we have 
\be
\partial_t \log |\varphi |^2 = {1 \over 6}\partial_t \log {\rm det} \, \tilde{g}
\ee
Then, we conclude
\be \label{evol-dilaton0}
\p_t \log | \varphi |^2 = | \varphi|^2 \big\{ -2 \nabla_\mu \alpha^\mu - R + 2 |\alpha|^2  \big\}
\ee
The flow of $g_{ij} = | \varphi |^{-2} \tilde{g}_{ij}$ is
\be
\partial_t g_{ij} = | \varphi |^{-2} \{ \partial_t \tilde{g}_{ij} - (\partial_t \log | \varphi |^2) g_{ij} \}.
\ee
Substituting the equations derived above,
\bea \label{evol-g}
\partial_t g_{ij} &=& - | \varphi|^2 \bigg\{ 2 ({\frak D}_p N_{ij}{}^p  +  {\frak D}_p N_{ji}{}^p) - \nabla_i \nabla_j \log | \varphi |^2 + J^p{}_i J^q{}_j \nabla_p \nabla_q \log | \varphi |^2 \nonumber\\
&&  - \alpha_i \alpha_j  + \alpha_{Ji} \alpha_{Jj}   +4 \alpha_p (N_j{}^p{}_i + N_i{}^p{}_j) \bigg\}
\eea
using $\alpha_i = - \partial_i \log | \varphi |^2$. The Ricci
curvature of $g_{ij}$ is given by (\ref{ricci-typeiia}).
Substituting this into (\ref{evol-g}), we obtain the flow of $g_{ij}$ as stated in Theorem \ref{main}. Q.E.D.

\bigskip

\noindent Department of Mathematics $\&$ Computer Science, Rutgers, Newark, NJ 07102, USA

\smallskip

\noindent teng.fei@rutgers.edu

\bigskip

\noindent Department of Mathematics, Columbia University, New York, NY 10027, USA

\smallskip

\noindent phong@math.columbia.edu

\bigskip

\noindent \noindent Mathematics Department, University of British Columbia, Vancouver, BC V6T 1Z2, CAN

\smallskip

\noindent spicard@math.ubc.ca

\bigskip

\noindent Department of Mathematics, University of California, Irvine, CA 92697, USA

\smallskip
\noindent xiangwen@math.uci.edu

\end{document}